\documentclass {svjour2}                    
\smartqed  
\usepackage{graphicx}
\newcommand{\eeq}{\end{equation}}
\newcommand{\beq}{\begin{equation}}
\newcommand{\nuq}[1]{\label{#1} \eeq}
\newcommand{\ba}{\begin{array}}
\newcommand{\ea}{\end{array}}
\newcommand{\veps}{\varepsilon}

\newcommand{\gtle}{\mbox{$\stackrel{>}{\scriptstyle <}$}}
\newcommand{\legt}{\mbox{$\stackrel{<}{\scriptstyle >}$}}
\journalname{Journal of Statistical Physics}

\begin{document}

\title{The global statistics of return times:
return time dimensions versus generalized measure dimensions}

\author{Giorgio Mantica}

\institute{G. Mantica \at Center for Nonlinear and Complex
Systems, Dipartimento di fisica e matematica, Universit\`a
dell'Insubria, Via Valleggio 11, 22100 Como, Italy and INFN,
sezione di Milano, and CNISM, Unit\`a di Como.
\email{giorgio@uninsubria.it} }

\maketitle

\begin{abstract}
We investigate the relations holding among generalized dimensions
of invariant measures in dynamical systems and similar quantities
defined by the scaling of global averages of powers of return
times. Because of a heuristic use of Kac theorem, these latter
have been used in place of the former in numerical and
experimental investigations; to mark this distinction, we call
them return time dimensions. We derive a full set of inequalities
linking measure and return time dimensions and we comment on their
optimality with the aid of two maps due to von Neumann -- Kakutani
and to Gaspard -- Wang. We conjecture the behavior of return time
dimensions in a typical system. We only assume ergodicity of the
dynamical system under investigation.


\end{abstract}

\section{Introduction: Statement of the Problem and Previous Results}
\label{secintro}

The metric theory of dynamical systems is based on the study of a
transformations $T$ of a space $X$ into itself, that preserves a
probability measures $\mu$ on a suitable sigma algebra $\mathcal A$. We
assume throughout this paper that the dynamical system $(X,T,{\mathcal
A},\mu)$ under investigation is ergodic and that $X$ is a compact metric
space enclosed in ${\bf R}^n$. This case is general enough to cover many
practical applications. From a physical point of view, $\mu$ may be
thought of as the invariant distribution in the space $X$ of points of a
typical trajectory of the system, generated by repeated applications of
the map $T$ on a starting point $x$.

Our interest lies in {\em return times} of the motion. Let $A \in
{\mathcal A}$ be a measurable subset of $X$ of positive measure.
Later on, we shall choose $A$ to be a ball, that is, a circular
neighborhood of a point. Let $x$ be any point in $A$. We denote by
$\tau_A(x)$ the (integer) time of the first return of $x$ to the
set $A$:
\begin{equation}
 \tau_A(x) = \inf \{ n>0 \; \mbox{s.t.} \;  T^n(x)\in A \} .
\label{rito1}
\end{equation}

Return times and invariant measures are linked by a variety of
results that stand on the pillars of the classical theorems of
Poincar\'e and Kac \cite{peters}. The first guarantees that the
return time of a point $x$ to the set $A$ is almost surely finite,
with respect to any invariant measure $\mu$; the second links the
average time needed for recurrence of points in the set $A$ to the
inverse of the measure of $A$. On these bases, it was conjectured
long ago by Grassberger \cite{grass} and independently by Jensen
{\em et al.} \cite{kada} that the statistical moments of return
times, when {averaged} over balls of radius $\veps$, centered at
{\em all} points of a typical trajectory (therefore, not uniquely
fixed as in Poincar\'e and Kac theorems), have a power--law
scaling behavior, when $\veps$ tends to zero, with exponents
proportional to the generalized dimensions of the measure $\mu$.

Generalized dimensions of measures, defined {\em \`a la}
Hentschel--Procaccia
\cite{hents,grassp,grass,sandro,remo1,remo2,cutler}, have a large
importance in dynamical systems, see Pesin \cite{pesin-book} for a
comprehensive review. Their computation is a task of practical and
theoretical relevance, for which many alternative techniques have
been proposed. Therefore, Grassberger and Jensen {\em et al.} idea
offers a most interesting alternative in this respect.

Indeed, the original conjecture of  has become implicit usage in
successive investigations, that have computed generalized
dimensions from the statistics of return times. Yet, even before
the most recent applications of this technique
\cite{nat-0,nature}, this approach has been critically examined in
\cite{sandold}. Stimulated by these findings, we have tried to
answer a fundamental question that has frequently been overlooked:
whether the conjecture is rigorous and whether it is exact in
certain cases, the former obviously implying the latter. In order
to disambiguate this point, in this paper we shall call {\em
return time dimensions} the values obtained from the scaling
exponent of averages of return times, and we shall investigate
whether they are equal to {\em measure generalized dimensions}.

Before getting into details, observe that the conjecture is bold:
generalized measure dimensions are defined independently of the
dynamics, while return times obviously are. Put in another way,
the same measure (characterized by a spectrum of
Hentschel--Procaccia--Pesin generalized dimensions) can be the
equilibrium measure of quite different dynamical systems.
Precisely because of this, studying the relations holding among
the two sets of dimension is interesting, independently of the
validity of the above conjecture, since it leads to ``universal''
results that hold for all dynamical maps $T$ for which a given
measure $\mu$ is invariant.

In a first paper \cite{prl} we have studied this problem for
invariant measures supported on attractors of Iterated Function
Systems. The scope of this work has been successively enlarged in
\cite{ugacha} by the analysis of return (and entrance) times in
dynamical cylinders (rather than balls) for Bowen--Gibbs measures.
Relying on precise approximations to the {\em local} statistics of
return times obtained in \cite{miguel} the situation for entrance
times (a variant of the approach mentioned above) has been almost
completely clarified, while that for return times has been settled
only for indices $q<1$ (see below for definitions and further
discussions). Results concerning return times in cylinders and
their fluctuations are numerous: see {\em e.g.} \cite{colgas},
\cite{beno}. For the class of super--disconnected I.F.S. cylinders
and balls are in a strict relation, described in \cite{prl}. Yet,
in the general case, the problem of return times in balls, rather
than cylinders, remains completely open and it is arguably the
most relevant to practical and numerical applications.

In this paper we advance the analysis of this problem by proving
rigorous bounds holding in full generality between measure
dimensions and those obtained via return times. In fact, we do not
require any additional property (like {\em e.g.} Bowen--Gibbs) on
the dynamical system under consideration, other than those listed
above. In the course of this analysis we will also consider the
comparison between generalized dimensions and their box versions,
commonly used in numerical simulations. We shall introduce new box
quantities which will be shown to be optimal, both for measure and
for return time dimensions. Finally, by analyzing the case of two
significant one--dimensional maps, we shall demonstrate the
optimality of the derived inequalities and we shall put forward a
conjecture on the behavior of return times dimensions in a
``typical'' case.

On the contrary, we shall not consider the problem of the
multifractal decomposition, {\em i.e.} whether dimensions are
linked to the so--called $f(\alpha)$ spectrum
\cite{multif1,pesjstat,olsen}.

It must finally be underlined the difference of this problem---the
{\em global} statistics of return times---with the much more
investigated case of the {\em local} statistics, that consider the
distribution of return times of points in a nested sequence of
neighborhoods of a given point: see {\em e.g.}
\cite{ornstein,beno,beno2,nico,miguel} and references therein.

\section{Definitions, Structure of the Paper and Summary of Results}
\label{secdefs} We start by giving formal definitions of
generalized dimensions (a variety of possibilities are encountered
in the literature). Let $B_\varepsilon(x)$ be the ball of radius
$\varepsilon$ at $x$ and $q$ a real quantity different from one.
The {\em partition functions} $\Gamma_\mu(\varepsilon,q)$ and
$\Gamma_\tau(\varepsilon,q)$ are the integrals
\begin{equation}
 \Gamma_\mu(\varepsilon,q) := \int_X [\mu({B_\varepsilon(x)})]^{q-1}
 d\mu(x),
 \label{nte2}
 \end{equation}
\begin{equation}
 \Gamma_\tau(\varepsilon,q) := \int_X [\tau_{B_\varepsilon(x)}(x)]^{1-q}
 d\mu(x).
 \label{nte1}
 \end{equation}
If the integrand is not summable, we shall understand that the
value of the partition function is infinite. The symmetry between
the two definitions is apparent and betrays the idea behind the
approach mentioned in the Introduction: the measure of a ball in
eq. (\ref{nte2}) is replaced in eq. (\ref{nte1}) by the inverse of
the return time of the point at its center. Remark that the
integral in eq. (\ref{nte1}) can be computed by a Birkhoff sum
over a trajectory \cite{sandold,prl}, as in the original proposals
\cite{grass,kada}. Remark also that the actual numerical
computations for \cite{grass} were performed with Birkhoff sums of
the kind $ \sum_{i,j} \tau_{B_\varepsilon(x_i)}^{1-q}(x_j)$ (P.
Grassberger, private communication) and therefore they were
estimates of the integral $\int [\tau_{B_\varepsilon(x)}(y)]^{1-q}
d\mu(x) d\mu(y)$. This amounts to computing entrance (rather than
return) times.

The {\em generalized dimensions} $D^\pm_\sigma(q)$ are defined via
the scaling of partition functions for small $\varepsilon$:
\begin{equation}
 \Gamma_\sigma(\varepsilon,q) \sim \varepsilon^{D^\pm_\sigma(q)(q-1)},
 \label{nte3}
 \end{equation}
where  $\sigma$ from now on denotes either $\mu$ or $\tau$. More
precisely, one has that
\begin{equation}
  D^\pm_\sigma(q) :=   \limsup (\inf)
  \frac{1}{q-1}
   \frac{\log \Gamma_\sigma(\varepsilon,q)}
  {\log \varepsilon}.
 \label{nte3b}
 \end{equation}
For $q=1$, as usual, slightly different definitions are needed:
\begin{equation}
  D^\pm_\sigma(1) :=   \limsup (\inf) \frac{\Gamma^l_\sigma(\varepsilon)}
  {\log \varepsilon},
 \label{nte3l1}
 \end{equation}
where
\begin{equation}
  \Gamma^l_\mu(\varepsilon) := \int_X \log [\mu({B_\varepsilon(x)})],
 \label{nte3l2}
 \end{equation}
\begin{equation}
  \Gamma^l_\tau(\varepsilon) :=   \int_X \log [\tau^{-1}_{B_\varepsilon(x)}(x)].
 \label{nte3l3}
 \end{equation}
As noted, partition functions may be infinite: in such case we
shall also set equal to infinity the corresponding generalized
dimensions.

The central question addressed in this paper is the nature of the
relations between the two sets of dimensions: are they equal,
always or in certain cases at least~? Can a set of rigorous
inequalities among them be derived~?

The results of this paper are organized as follows. In the next
section we briefly outline basic properties (monotonicity,
convexity) of return time generalized dimensions. Then, we shall
find it convenient to introduce a number of additional quantities,
that we shall also call dimensions and that are interesting on
their own. Some of these dimensions are conventional, some are
new. In Sect. \ref{secboxdims} we start by defining box
dimensions, both for measure and return times, following typical
usage in experimental and numerical applications: a {\em grid} of
box-partitions of the space is considered, and limits are taken
with respect to this grid.

In Sect. \ref{sectriedi} we review the known relations between box
and generalized dimensions of measures and a proposal put forward
by Riedi \cite{rie} to avoid ``pathological'' values of box
dimensions for negative $q$. By a modification of his idea we
define a new box partition function that offers a definite
theoretical advantage over both the original box quantities and
Riedi's enhanced box formalism: its scaling yields the generalized
dimensions $D_\mu^\pm(q)$ for all values of $q$, independently of
the particular grid adopted. This is made formal in Theorem
\ref{teodima}, that, although not directly related to return
times, constitutes one of the main results of this paper.

While the previous results deal with measure dimensions, in Sect.
\ref{bvcrtdim} we consider the relations between generalized and
box dimensions for return times: Proposition \ref{prop2} shows
that the latter are always larger than, or equal to, the former.
Mimicking the procedure developed for measures, we introduce a box
quantity that yields exactly the generalized return time
dimensions $D^\pm_\tau(q)$, again for all values of $q$,
independently of the particular grid adopted: this is the content
of Theorem \ref{teodimab}.

We then put in relation measure and return time dimensions,
according to the theme of this paper. Section \ref{secdiscret}
introduces a central quantity to this goal: the distribution of
return times into a fixed set $A$ in the space $X$. The zeroth and
first moment of this distribution are fixed by Poincar\'e and Kac
theorems. Basic inequalities are derived for the remaining
moments: Lemma \ref{lema1}. We stress again that this is obtained
in the most general setting.

These results are put at work in Sect. \ref{secrelbox}:
inequalities between measure and return box dimensions are derived
for all values of $q$ and equality is found for $q=0$: Proposition
\ref{prop1} gives full detail.

Section \ref{secallthings} is the heart of the paper. Here, we
chain together our results in Theorem \ref{teototu}, that presents
the most complete set of inequalities holding in full generality
among generalized and box dimensions, for measures and return
times. In the same section, we discuss the optimality of the
inequalities presented. We study the role of short returns, that
imply an upper bound for positive dimensions, Lemma
\ref{lemdeca1}. We also investigate the different situations
occurring for positive and negative values of $q$ and we compute
return time dimensions in two interesting cases: in full detail
for the von Neumann -- Kakutani map \cite{von} (Theorem
\ref{teokaku}) and, partly, for the Gaspard -- Wang intermittent
map \cite{gasp} (Theorem \ref{teomane}). Also in this section we
formulate a conjecture on the typical behavior of return time
generalized dimensions that links it significantly to measure
generalized dimensions.

Conclusions are presented briefly in Sect. \ref{secconcl}, while
three additional sections, \ref{secmapintex}, \ref{secprteokaku}
and \ref{secintermap} contain the details of the calculations and
proofs for the two maps quoted above, as well as side results.

\section{General properties of return time dimensions}
\label{secgenprop}

Because of the formal similarity between eqs.  (\ref{nte2}) and
(\ref{nte1}) some of the properties of generalized measure
dimensions also characterize return time dimensions. A couple of
these are contained in the following Lemma.

\begin{lemma}
The return time generalized dimensions $D^{\pm}_\tau(q)$ are
monotone non increasing in the index $q$ and the functions $(q-1)
D^{-}_\tau(q)$ for $q>1$ and $(q-1) D^{+}_\tau(q)$ for $q<1$ are
convex. \label{lemgprop}
\end{lemma}
{\em Proof.} Observe that both $\Gamma_\mu(\varepsilon,q)$ and
$\Gamma_\tau(\varepsilon,q)$ can be seen as integral of a function
$\phi(x)$ raised to the power $q-1$. In the return time case this
latter is $\phi(x) = 1/\tau_{B_\varepsilon(x)}(x)$. The two
results above are then a consequence of Jensen and Holder
inequalities, similar to those holding for $D^{\pm}_\mu(q)$, whose
details can be found in \cite{belli,tre}. \qed

Additional results can be obtained in this line, but will be
reported elsewhere. In fact, our specific aim in this paper is
simply to compare the value of measure and return time dimensions.
In this regard, finiteness of the return time dimensions is an
important issue that will be considered in Sections \ref{bvcrtdim}
and \ref{secallthings}.

\section{Box Dimensions of Measures and of Return Times Distributions}
\label{secboxdims} Usually, in numerical experiments, rather than
computing the integral (\ref{nte2}) one covers the set $X \subset
{\bf R}^n$ by a lattice of hypercubic boxes $A_j$, $j=0,1,\ldots$
of side $\varepsilon$. The usual choice is to draw the zeroth box
as having the origin of the coordinates as a corner and the sides
exiting from that corner oriented as the coordinate directions.
Clearly, different choices are possible, varying origin and
orientation. We shall let $\theta \in \Theta$ denote the
particular choice of origin and orientation in the set $\Theta$ of
all choices. $\theta$ will be called a {\em grid}. Therefore, a
grid consists of infinitely many box partitions of $X$, one for
every value of $\veps$.

We then consider in place of the partition functions
$\Gamma_\mu(\varepsilon,q)$, the sums
\begin{equation}
 \Upsilon_\mu(\theta,\varepsilon,q) :=
  \sum_{j  \; \mbox{\scriptsize s.t.} \; \mu({A_j}) > 0}
  \mu({A_j})^{q}.
 \label{nte2b}
 \end{equation}
For simplicity of notation, the dependence of $A_j$ on $\theta$
and $\veps$ will be left implicit here and in the following.

Similarly, by replacing in eq. (\ref{nte1}) the box centered at
$x$ with the set $A_j$ that contains $x$, we define the {\em
return time box partition function}:
\begin{equation}
 \Upsilon_\tau(\theta,\varepsilon,q) := \sum_j
  \int_{A_j} \tau_{A_j}^{1-q}(x) d\mu(x).
 \label{nte1b}
 \end{equation}
It has to be noticed the double role of the set $A_j$, as starting
and arrival set of the motion. We shall break this symmetry later
on. The logarithmic analogues are
\begin{equation}
 \Upsilon^l_\mu(\theta,\varepsilon) := \sum_{j  \; \mbox{\scriptsize s.t.} \; \mu({A_j}) > 0}
  \mu({A_j}) \log [\mu({A_j})],
 \label{nte2x}
 \end{equation}
and
\begin{equation}
 \Upsilon^l_\tau(\theta,\varepsilon) :=   \sum_j
  \int_{A_j} \log [\tau^{-1}_{A_j}(x)] d\mu(x).
 \label{nte1x}
 \end{equation}

Define now the {\em box generalized dimensions}
$\Delta^\pm_\mu(\theta,q)$ and $\Delta^\pm_\tau(\theta,q)$, by
using $\Upsilon$'s and $\Delta$'s in place of $\Gamma$'s and
$D$'s, respectively, in eqs. (\ref{nte3b}) and (\ref{nte3l1}):

\begin{equation}
  \Delta^\pm_\sigma(\theta,q) :=   \limsup (\inf)
  \frac{1}{q-1}
   \frac{\log \Upsilon_\sigma(\theta,\varepsilon,q)}
  {\log \varepsilon},
 \label{nte3brt}
 \end{equation}
\begin{equation}
  \Delta^\pm_\sigma(\theta,1) :=   \limsup (\inf) \frac{\Upsilon^l_\sigma(\theta,\varepsilon)}
  {\log \varepsilon},
 \label{nte3l1rt}
 \end{equation}
where $\sigma$ can be either $\mu$ or $\tau$.

Notice that it is possible to avoid dependence on the specific
grid by taking the {\em infimum} (or the supremum, according to
the value of $q$) over all such grids in the definition of the
partition function $\Upsilon_\tau(\varepsilon,q)$
\cite{pesin-book}. We elect not to take this step for two reasons.
The first is that this is difficultly achievable in numerical
applications and for the same reason it cannot be a good model of
what is numerically observed. The second is that we shall strive
at obtaining results that do not depend on the particular grid
selected, but apply to all and {\em a fortiori} also to dimensions
defined with the infimum procedure included. A first instance of
this fact is to be met in the next section, where we introduce
enhanced box dimensions, based on an idea pioneered by Riedi.

\section{Box versus Generalized Measure Dimensions}
\label{sectriedi} The relations among different measure dimensions
is a subject that has been intensively studied, see {\em e.g.}
\cite{pesjstat},\cite{guy},\cite{pesin-book},\cite{rie},\cite{tre},\cite{due}
with an effort towards proving their equivalence, on one side, and
towards releasing the request of performing an infimum over grids,
on the other side. We first review the known relations needed for
our scope in the following Lemma, and then we present a new result
that we believe to be of some importance.
\begin{lemma}
The following relations exists between box and generalized measure
dimensions:
 \beq
  \ba{ll}
 \Delta^{\pm}_\mu(\theta,q)  =   D^{\pm}_\mu(q)   & \mbox{ for }  q > 0, \\
 \Delta^{\pm}_\mu(\theta,q)  \geq   D^{\pm}_\mu(q)   & \mbox{ for }  q \leq 0. \\
 \ea
\nuq{ups8x} \label{lemaboxc1}
\end{lemma}
{\em Proof.} The full proof, including the non-trivial interval $q
\in (0,1]$, can be found in the complete exposition \cite{tre}.
Notice that no infimum procedure over $\veps$-grids is involved.

Therefore, measure box dimensions are independent of the choice of
the grid $\theta$ for any $q>0$. Examples exist showing both such
dependence and strict inequality w.r.t. generalized dimensions for
$q<0$. The case $q=0$ seems to be on less firm ground, see Sect.
\ref{secallthings}. Roughly speaking, what might happen for
negative $q$ is the following: if a box $A_j$ ``barely touches''
the support of the measure $\mu$ close to one of its edges, its
measure can be arbitrarily small, independently of its size
$\veps$, so that $\Delta^{+}_\mu(\theta,q)$ can become arbitrarily
large. To avoid this effect Riedi \cite{rie} introduced the sums
\begin{equation}
 \Phi_\mu(\theta,\varepsilon,q) := \sum_{j  \; \mbox{\scriptsize s.t.} \; \mu({A_j}) > 0}
  \mu({\overline{A}_j})^{q},
 \label{nte4b}
 \end{equation}
where $\overline{A}_j$ is a box of side $3 \varepsilon$ centered
on the box $A_j$, cfr. eq. (\ref{nte2b}). In one dimension, for
instance, $\overline{A}_j$ consists of the union of $A_{j-1}$,
$A_j$ and $A_{j+1}$. The geometrical situation in two and more
dimensions can be easily pictured by the reader.

Using clever manipulations, Riedi has been able to prove that the
generalized dimensions generated by the scaling of $\Phi_\mu$
coincide with $D_\mu^{\pm}(q)$ for $q>1$. It actually follows from
estimates in \cite{tre} that equality can be proven for any $q>0$.
Although it is plausible that this also holds in large generality
for negative $q$ as well (see the numerical results in
\cite{sator}), we have not been able to find a formal proof of
this fact, that would hold in the most general setting adopted in
this paper. Yet, in this endeavor, we have discovered a new box
quantity that achieves this goal:

\begin{theorem} For any $\theta \in \Theta$,
the scaling behavior of the partition function defined via:
\begin{equation}
 \Psi_\mu(\theta,\varepsilon,q) := \sum_{j  \; \mbox{\scriptsize s.t.} \; \mu({A_j}) > 0}
  \mu({A_j}) \mu({\overline{A}_j})^{q-1},
 \label{nte4c}
 \end{equation}
for $q \neq 1$ and
\begin{equation}
 \Psi^l_\mu(\theta,\varepsilon) := \sum_{j  \; \mbox{\scriptsize s.t.} \; \mu({A_j}) > 0}
  \mu({A_j}) \log [\mu({\overline{A}_j})],
 \label{nte4c2}
 \end{equation}
for $q=1$ yields the generalized dimensions $D_\mu^{\pm}(q)$.
\label{teodima}
\end{theorem}

Proof of Thm. \ref{teodima} is a direct consequence of the
auxiliary results collected in:
\begin{lemma}
For any $\theta \in \Theta$ the following inequalities hold:
  \beq
  \Gamma_\mu(\varepsilon,q) \leq \Psi_\mu(\theta,\varepsilon,q)
  \leq \Phi_\mu(\theta,\varepsilon,q),
  \;\;\; q \geq 1
 \nuq{riedii}
  \beq
  \Gamma_\mu(\varepsilon,q) \geq \Psi_\mu(\theta,\varepsilon,q)
  \leq \Phi_\mu(\theta,\varepsilon,q),
  \;\;\; q \leq 1
 \nuq{riedi2}
\beq
  \Gamma_\mu(k \varepsilon,q) \leq \Psi_\mu(\theta,\varepsilon,q)
  \;\;\; q \leq 1,
 \nuq{mio2}
 \beq
  \Gamma_\mu(k \varepsilon,q) \geq \Psi_\mu(\theta,\varepsilon,q)
  \;\;\; q \geq 1,
 \nuq{mio2b}
\beq
  \Gamma^l_\mu(k \varepsilon) \geq \Psi^l_\mu(\theta,\varepsilon)
   \geq \Gamma^l_\mu(\varepsilon),
 \nuq{mio2c}

where $k$ is a dimension-dependent multiplier. \label{lema22}
\end{lemma}

{\em Proof.} To prove the first inequality we follow \cite{rie}.
Since
\begin{equation}
\Gamma_\mu(\varepsilon,q) := \int_X
d\mu(x)[\mu({B_\varepsilon(x)})]^{q-1} = \sum_j \int_{A_j} d\mu(x)
[\mu({B_\varepsilon(x)})]^{q-1}
 \label{mio3}
 \end{equation}
and since $B_\varepsilon(x) \subset \overline{A}_j$ when $x \in
A_j$, $q \geq 1$,
\begin{equation}
\Gamma_\mu(\varepsilon,q) \leq \sum_{j } \int_{A_j} d\mu(x)
[\mu(\overline{A}_j)]^{q-1} = \Psi_\mu(\theta,\varepsilon,q) \leq
\Phi_\mu(\theta,\varepsilon,q).
 \label{mio4}
 \end{equation}
Equally,
\begin{equation}
\Gamma^l_\mu(\varepsilon) \leq \sum_{j } \int_{A_j} d\mu(x) \log
[\mu(\overline{A}_j)] = \Psi^l_\mu(\theta,\varepsilon).
 \label{mio4b}
 \end{equation}

Conversely, when $q \leq 1$, the first inequality in (\ref{mio4})
is reversed, while the second still holds:
\begin{equation}
\Gamma_\mu(\varepsilon,q) \geq \sum_{j  } \int_{A_j} d\mu(x)
[\mu(\overline{A}_j)]^{q-1} = \Psi_\mu(\theta,\varepsilon,q) \leq
\Phi_\mu(\theta,\varepsilon,q).
 \label{mio4b2}
 \end{equation}
Next, we use the fact that $\overline{A}_j \subset B_{k
\varepsilon}(x)$ when $x \in A_j$ and $k$ is a fixed multiplier,
as in section \ref{bvcrtdim}, to obtain, still for $q \leq 1$
\begin{eqnarray}
\Gamma_\mu(k \varepsilon,q) = \sum_{j } \int_{A_j} [\mu(B_{k
\varepsilon}(x) ]^{q-1}d\mu(x) \leq \nonumber \\
 \leq \sum_{j } \int_{A_j}
[\mu(\overline{A}_j)]^{q-1}d\mu(x) =
 \sum_{j  } \mu({A}_j)
[\mu(\overline{A}_j)]^{q-1}
  = \Psi_\mu(\theta,\varepsilon,q).
 \label{mio5}
 \end{eqnarray}
Finally, for $q \geq 1$, we get
\begin{eqnarray}
\Gamma_\mu(k \varepsilon,q)  \geq \sum_{j } \int_{A_j}
[\mu(\overline{A}_j)]^{q-1}d\mu(x) =
\Psi_\mu(\theta,\varepsilon,q),
 \label{mio5x}
 \end{eqnarray}
\begin{equation}
\Gamma_\mu^l(k \varepsilon) \geq \sum_{j  } \int_{A_j} d\mu(x)
\log [\mu(\overline{A}_j)] = \Psi^l_\mu(\theta,\varepsilon).
 \label{mio4x} \qed
 \end{equation}

Theorem \ref{teodima} asserts that the grid--dependent sums
$\Psi_\mu(\theta,\varepsilon,q)$ give rise to a set of dimensions
that are independent of the grid $\theta$ and coincide with the
generalized dimensions for all values of $q$. In addition, from a
numerical point of view, the sums $\Psi_\mu(\theta,\varepsilon,q)$
can be evaluated with the same effort required for computing the
original sums  (\ref{nte2b}) or Riedi's extension (\ref{nte4b}):
all one needs to know is the value of the box measures $\mu(A_j)$.
We therefore believe that Theorem \ref{teodima} can become a new
tool in the multifractal analysis of measures.

\section{Box versus Generalized Return Time Dimensions}
\label{bvcrtdim} In a completely analog way to what done in the
previous section for measure dimensions, we now compare the box
and generalized return time dimensions,
$\Delta^{\pm}_\tau(\theta,q)$ and $D^{\pm}_\tau(q)$. Then, we
introduce a new box quantity for return times, $\Psi_\tau$,
analogous to $\Psi_\mu$ of the previous section.  We prove that
this box partition function yields the generalized return time
dimensions $D^{\pm}_\tau(q)$ for all values of $q$, independently
of the grid $\theta$.

\begin{proposition}
The box return time dimensions $\Delta^{\pm}_\tau(q)$ are always
larger than, or equal to, their generalized counterparts:
  \beq
 \Delta^{\pm}_\tau(\theta,q)  \geq   D^{\pm}_\tau(q).
\nuq{ups7}
\label{prop2}
\end{proposition}
{\em Proof.} Fix a specific grid $\theta$ and let $j(x)$ be the
index of the the hypercube of side $\varepsilon$ containing the
point $x$. Then, $A_{j(x)}$ is enclosed in the ball of radius $k
\varepsilon$ centered at $x$, with a fixed multiplier $k \geq 1$
that can be chosen as a function only of the (Euclidean) dimension
of the space. This implies that $\tau_{B_{k \varepsilon}(x)} (x)
\leq \tau_{A_j(x)}(x)$ for all $x$. Therefore, we part the
integral defining $\Gamma_\tau(k \varepsilon,q)$ over the grid of
side $\varepsilon$,
\begin{equation}
 \Gamma_\tau(k \varepsilon,q) = \int_X \tau_{B_{k\varepsilon}(x)}^{1-q}(x)
 d\mu(x) = \sum_j \int_{A_j} \tau_{B_{k\varepsilon}(x)}^{1-q}(x)
 d\mu(x),
 \label{nte1x1}
 \end{equation}
  and we use this inequality, first for $1-q \geq 0$, to get
 \begin{equation}
  \Gamma_\tau(k\varepsilon,q) \leq \sum_j \int_{A_j} \tau_{A_j}^{1-q}(x)
 d\mu(x) = \Upsilon_\tau(\theta,\varepsilon,q).
 \label{nte1x2}
 \end{equation}
For $1-q \leq 0$ we obtain the reverse inequality. In force of
these inequalities, an immediate calculation provides the thesis.
As before, the case $q=1$ requires a separate treatment:
\begin{equation}
  \Gamma^l_\tau(k \varepsilon)  =
   \sum_j \int_{A_j} \log [\tau^{-1}_{B_{k \varepsilon(x)}}(x)] d\mu(x)
  \geq    \sum_j \int_{A_j} \log [\tau^{-1}_{A_j}(x)] d\mu(x) =
 \Upsilon^l_\tau(\theta,\varepsilon).
 \label{nte42}
\end{equation}
Using this information in the limits (\ref{nte3b}),(\ref{nte3l1})
yields the thesis.
 \qed

The estimate in the previous proof help us also to establish
existence of the return time partition functions for $q \geq 0$.

\begin{lemma}
The the partition sums $\Gamma_\tau(\varepsilon,q)$ and
$\Upsilon_\tau(\theta,\varepsilon,q)$, as well as
$\Upsilon^l_\tau(\theta,\varepsilon)$ and
$\Gamma^l_\tau(\varepsilon)$ exist for any $\theta \in \Theta$, $q
\geq 0$.
 \label{lem-existence}
\end{lemma}
{\em Proof.} Existence is trivial for $q > 1$. If $0 \leq q \leq
1$ the functions $\Upsilon_\tau(\theta,\varepsilon,q)$ and
$\Upsilon^l_\tau(\theta,\varepsilon)$ exist because of Kac theorem
\cite{kac}. Then, the inequality (\ref{nte1x2}), valid for $q \leq
1$ and the inequality (\ref{nte42}) imply that also
$\Gamma_\tau(\varepsilon,q)$ and $\Gamma^l_\tau(\varepsilon)$
exist. \qed

The ideas exploited in the previous section can also be used to
construct a box quantity capable of generating the generalized
return time dimensions $D_\tau(q)$. This is defined via
\begin{equation}
 \Psi_\tau(\theta,\varepsilon,q) := \sum_j
  \int_{A_j} \tau_{\overline{A}_j}^{1-q}(x) \; d\mu(x),
 \label{nted3}
 \end{equation}
and
\begin{equation}
 \Psi^l_\tau(\theta,\varepsilon) := \sum_j
  \int_{A_j} \log(\tau^{-1}_{\overline{A}_j}(x)) \; d\mu(x).
 \label{nted3l}
 \end{equation}
Difference with eq. (\ref{nte1b}) has to be appreciated: the
integral is taken over the set $A_j$, but the return time is
computed when $x$ gets back into the larger set $\overline{A}_j$,
defined as in Sect. \ref{sectriedi}.

\begin{theorem} The scaling behavior of the
partition functions $\Psi_\tau(\theta,\varepsilon,q)$ and
$\Psi^l_\tau(\theta,\varepsilon)$ yield the generalized dimensions
$D_\tau^{\pm}(q)$ for any $\theta \in \Theta$. \label{teodimab}
\end{theorem}

{\em Proof.} Let again $j(x)$ be the index of the the hypercube of side
$\varepsilon$ containing the point $x$. The key point is that
\begin{equation}
   B_\varepsilon (x) \subset  \overline{A}_{j(x)}  \subset B_{ k\varepsilon} (x)
 \label{nted4}
 \end{equation}
with a dimension-dependent constant $k$. Therefore,
\begin{equation}
   \tau_{B_\varepsilon (x)}(x) \geq
   \tau_{\overline{A}_{j(x)}}(x)   \geq \tau_{ B_{ k\varepsilon}
   (x)}(x),
 \label{nted5}
 \end{equation}
which leads to
\begin{equation}
 \Gamma_\tau(\varepsilon,q) =  \sum_j \int_{A_j} \tau_{B_{\varepsilon}(x)}^{1-q}(x)
 d\mu(x) \leq  \Psi_\tau(\theta,\varepsilon,q) \leq \Gamma_\tau(k \varepsilon,q)
 \label{nted6}
 \end{equation}
for $q \geq 1$ and to a reverse chain of inequalities when $q \leq
1$. The logarithmic partition function, to be used for $q=1$
satisfies
\begin{equation}
 \Gamma^l_\tau(\varepsilon) =  \sum_j \int_{A_j} \log
 (\tau^{-1}_{B_{\varepsilon}(x)}(x))
 d\mu(x) \leq
\sum_j \int_{A_j} \log
 (\tau^{-1}_{\overline{A}_{j(x)}}(x)) d\mu(x) =
 \Psi^l_\tau(\theta,\varepsilon),
 \label{ntedl}
 \end{equation}
 and
\begin{equation}
\Psi^l_\tau(\theta,\varepsilon)  \leq \sum_j \int_{A_j} \log
 (\tau^{-1}_{B_{k \varepsilon}(x)}(x)) d\mu(x) = \Gamma^l_\tau(k
 \varepsilon),
 \label{ntedlb}
 \end{equation}
from which the thesis follows. \qed

 Remark that the geometric relation in eq. (\ref{nted4})
shows that return time dimensions are somehow performing the same
kind of action implied in Riedi's enlarged box idea. Also remark
that one is free to chose a different grid $\theta$ at each value
of $\veps$.

We end this section by showing the existence of a particular
combination of $\Upsilon_\tau(\theta,\varepsilon,q)$ that also
yields the generalized dimensions for $q$ larger than one. This is
defined as follows. Fix a grid $\theta_0$ in ${\bf R}^d$ and a
value $\veps>0$. Let $e_i$, $i=1,\ldots,d$ unit orthogonal vectors
giving the direction of the grid. On this basis, construct $3^d$
parallel grids with the same directions of $\theta_0$: call them
$\theta_l$, $l=0,\ldots,3^d-1$. The first of these grids is the
original $\theta_0$, the others have origins shifted by lattice
vectors of the kind $\veps \sum_{i=1}^d n_i e_i$, where the $n_i$
can take the values $0,1,2$. For each of these grids, consider
boxes of side $3 \veps$, and on this basis, construct the box
partition function
 \beq
 \tilde \Upsilon_\tau(\theta_0,\varepsilon,q) :=
 \sum_{l=0}^{3^d-1}
 \Upsilon_\tau(\theta_l,3 \varepsilon,q).
 \label{utild1}
 \end{equation}

\begin{proposition}
When $q>1$, the box partition function $\tilde
\Upsilon_\tau(\theta_0,3 \varepsilon,q)$ yields the generalized
dimensions $D_\tau^\pm(q)$, independently of the choice of the
grid.
 \label{prop3}
\end{proposition}
{\em Proof.} Let $\theta_0$ be a given grid. Obviously, one has
\begin{equation}
\Psi_\tau(\theta_0,\varepsilon,q) := \sum_j
  \int_{A_j} \tau_{\overline{A}_j}^{1-q}(x) \; d\mu(x) \leq
\sum_j \int_{\overline{A}_j} \tau_{\overline{A}_j}^{1-q}(x) \;
d\mu(x),
 \label{nted7}
 \end{equation}
 where each integral has been extended to a larger domain.
The summation index $j$ runs over all boxes of size $\varepsilon$,
while the enlarged boxes $\overline{A}_j$ have side $ 3
\varepsilon$ and each of these is composed of $3^d$ smaller ones,
$d$ being the euclidean space dimension. Neighboring boxes
$\overline{A}_j$ overlap, but at the same time one can part the
$j$ summation into $3^d$ different sets of non-overlapping,
adjacent boxes. These are precisely defined by the $\theta_l$
grids defined above, so that eq. (\ref{nted7}) becomes
\begin{equation}
\Psi_\tau(\theta_0,\varepsilon,q)   \leq  \sum_{l=0}^{3^d-1}
\Upsilon_\tau(\theta_l,3 \varepsilon,q) = \tilde
\Upsilon_\tau(\theta_0, \varepsilon,q).
 \label{nted8}
 \end{equation}
The above equation is valid for all values of $q$. Let now $q
> 1$. Then, since $\Upsilon_\tau(\theta,\varepsilon,q) \leq
\Gamma_\tau(k\varepsilon,q)$ for any $\theta$ (see eq.
\ref{nte1x2}) and using also eq. (\ref{nted6}), we find
\begin{equation}
\Gamma_\tau(\veps,q) \leq \Psi_\tau(\theta_0,\varepsilon,q)   \leq
\sum_{l=0}^{3^d-1} \Upsilon_\tau(\theta_l,3 \varepsilon,q) =
\tilde \Upsilon_\tau(\theta_0,\varepsilon,q) \leq 3^d
\Gamma_\tau(3 k \veps,q).
 \label{ntre3}
 \end{equation}
The by-now usual technique proves the thesis. \qed

\section{Kac Theorem and moment inequalities}
\label{secdiscret} We need to bridge the gap between measure and
return time dimensions. To do this, our main tool will be Kac
theorem \cite{kac}, that we put at work in this section. For any
measurable set $A$ of positive measure, define the discrete return
times measure $\nu^A$ via
 \beq
 \nu^{A}(\{j\}) := \mu (
\{ x \in A \; \mbox{ s.t.} \;  \tau_A(x) = j \}) / \mu(A).
  \nuq{discrm1}
In words, $\nu^{A}(\{j\})$ is the normalized measure of the set of
points of $A$ that return to $A$ in $j$ time steps. Obviously,
$\nu^{A}$ is a measure supported on the positive integers, a fact
that will be exploited momentarily. Poincar\'e Theorem guarantees
that $\nu^{A}$ is a probability measure:
  \beq
 \sum_{j=1}^\infty \nu^{A}(\{j\}) = 1.
 \nuq{discrm0}
 We shall study the moments of this measure: for $s \in {\bf R}$,
let $\nu_s^A$ be:
  \beq
 \nu_s^A := \sum_{j=1}^\infty j^s \nu^{A}(\{j\}) =
 \frac{1}{\mu(A)}
  \int_A [\tau_{A}(x)]^{s} d\mu(x).
 \nuq{discrm2}
 Define also the logarithmic moment:
 \beq
 \nu_l^A := \sum_{j=2}^\infty \log (j)  \nu^{A}(\{j\}).
 \nuq{discrm4}

 Under the ergodicity hypothesis that we are assuming throughout,
Kac theorem fixes the value of the first moment of this measure:
 \beq
 \nu_1^A = 1 / \mu(A).
 \nuq{discrm3}
The key ingredient of our theory is the fact that all moments
$\nu_s^A$ can be put in relation to the latter, that to say, to
$\mu(A)$. In fact, we have the following Lemma.
\begin{lemma}
Let $\nu$ be a probability measure supported on $[1,\infty)$. Let
$\nu_s$ be its moments, allowing for an infinite value
of these latter. As a  function of $s$, $\nu_s$ is
monotonic, non--decreasing. Furthermore,
 \beq
  \ba{lll}
 \nu_s   \leq & (\nu_1)^s & \mbox{ for } 0 \leq s \leq 1 \\
 \nu_s   \geq & (\nu_1)^s & \mbox{ for } s \leq 0 , s \geq 1.
 \ea
\nuq{hold2} \label{lema1}
\end{lemma}

{\em Proof.} Since  $\nu$ is supported on $[1,\infty)$, monotonicity
follows immediately.  Apply H\"older inequality to $\nu_r$, $r \in {\bf
R}$, using the fact that $\nu$ is a probability measure:
 \beq
 |\nu_r| = |\int d\nu(x) x^r| \leq \| x^r \|_p \| 1 \|_{p'}
  = (\nu_{pr})^{\frac{1}{p}},
 \nuq{hold1x}
for any $r \in {\bf R}$ and any $p,{p'}>1$ such that
$p^{-1}+{p'}^{-1}=1$. Obviously, not all positive moments of
$\nu_r$ are finite, while this is true for all negative moments,
owing to the fact that $\nu$ is supported on $[1,\infty)$. This
also implies that all moments are positive, so that (\ref{hold1x})
can be simplified, to the extent that
  \beq
 (\nu_r)^p   \leq \nu_{pr}
 \nuq{hold1}
for any $r \in {\bf R}$ and any $p>1$.

 When $s$ is equal to
either zero, or one, equality of $\nu_s$ and $(\nu_1)^s$ holds
trivially. When $s \in (0,1)$, set $p=1/s$ and $r=s$ in
(\ref{hold1}), to get the thesis. When $s>1$, still use
(\ref{hold1}) letting $r=1$ and $p=s$. For negative values of $s$,
we start from the inequality:
 \beq
1 =  \int d\nu(x) = \int x^r x^{-r} d\nu(x)  \leq \| x^r \|_p \|
x^{-r} \|_{p'}
  = (\nu_{pr})^{\frac{1}{p}}
  (\nu_{-{p'}r})^{\frac{1}{p'}}.
 \nuq{hold3}
Letting $p=p'=2$ and $s=2r$ yields
 \beq
   \nu_{-s} \geq 1/ \nu_{s},
   \nuq{hold4}
valid for all real values of $s$.  We can now tackle the case
$s<-1$. Put $r=-1$ and $p=-s$ in (\ref{hold1}), to get
 \beq
   \nu_{s} \geq (\nu_{-1})^{-s} \geq (\nu_{1})^{s},
 \nuq{hold5}
 where the last inequality follows from $\nu_{-1} \geq
1/ \nu_{1}$, a particular case of (\ref{hold4}). Finally, for $s
\in (-1,0)$, use again (\ref{hold4}) and the first part of
(\ref{hold2}), that we have proven above and that applies since
$-s \in (0,1)$:
  \beq
   \nu_{s} \geq 1/ \nu_{-s}  \geq (\nu_{1})^{s}.
   \nuq{hold6}
This completes the proof. \qed

\begin{lemma}
In the same hypotheses of lemma \ref{lema1}, one has
 $ \nu_l:= \int \log(x) d\nu(x)   \leq   \log (\nu_1)$.
\label{lema2}
\end{lemma}
{\em Proof.} Since $\nu$ is a probability measure, this is Jensen's
inequality. \qed

Because of the observations made at the beginning of this section,
the above lemmas apply to $\nu_s^A$, the moments of the return
times of points in any positive measure set $A$, when taken with
respect to the normalized measure $d\mu_A(x) = \frac{1}{\mu(A)}
d\mu(x)$. As such, the formulae (\ref{hold2}) extend the content
of Kac theorem to all moments. Later in the paper, we shall find
examples where inequalities (\ref{hold2}) are strict, as well as
examples where they hold as equalities. We shall now investigate
the mathematical implications of these results to the dimension
problem.

\section{Inequalities between Measure and Return Times Box Dimensions}
\label{secrelbox}

On the basis of the theory of the previous section,  Lemmas
\ref{lema1} and \ref{lema2}, we can now study the quantities
$\Upsilon_\sigma(\theta,\varepsilon,q)$ and the associated
dimensions $\Delta^\pm_\sigma(\theta,q)$.
\begin{proposition}
The box dimensions $\Delta^{\pm}_\sigma(\theta,q)$,
$\sigma=\mu,\tau$, for any $\theta \in \Theta$ are linked by the
inequalities
 \beq
 \ba{ll}
 \Delta^{\pm}_\tau(\theta,q)  \geq \Delta^{\pm}_\mu(\theta,q)  & \mbox{ for }  q < 0, \\
\Delta^{\pm}_\tau(\theta,q)  \leq  \Delta^{\pm}_\mu(\theta,q)  & \mbox{ for }  q > 0, \\
\Delta^{\pm}_\tau(\theta,0)  =  \Delta^{\pm}_\mu(\theta,0).  &
 \ea
\nuq{ups5} \label{prop1}
\end{proposition}

{\em Proof.} Observe that
\begin{equation}
  \int_{A_j} \tau_{A_j}^{1-q}(x) d\mu(x) = \mu({A_j})
  \nu_{1-q}^{A_j}
 \label{ups2}
 \end{equation}
so that
\begin{equation}
 \Upsilon_\tau(\theta,\varepsilon,q) = \sum_{j  \;
  \mbox{\scriptsize s.t.} \; \mu({A_j}) > 0}
 \mu({A_j})  \nu_{1-q}^{A_j}.
 \label{ups3}
 \end{equation}
Therefore, using Lemma \ref{lema1} and Kac theorem, eq.
(\ref{discrm3}), we get
 \begin{equation}
 \Upsilon_\tau(\theta,\varepsilon,q)  \leq \sum_{j  \; \mbox{s.t.} \; \mu({A_j}) > 0}
 \mu({A_j})^q = \Upsilon_\mu(\theta,\varepsilon,q)
 \label{ups4}
 \end{equation}
for $q \in (0,1)$ and $\Upsilon_\tau(\theta,\varepsilon,q) \geq
\Upsilon_\mu(\theta,\varepsilon,q)$
in the opposite case. Using now eqs. (\ref{nte3},\ref{nte2}) we
can prove the two inequalities in (\ref{ups5}), for $q \neq 1$.
This latter can be treated by writing
 \begin{equation}
 \Upsilon^l_\tau(\theta,\varepsilon) = - \sum_{j  \;
 \mbox{\scriptsize s.t.} \; \mu({A_j}) > 0}
 \mu({A_j})  \nu_l^{A_j}.
 \label{ups39}
 \end{equation}
Using Lemma \ref{lema2}, we arrive at
$\Upsilon^l_\tau(\theta,\varepsilon) \geq
\Upsilon^l_\mu(\theta,\varepsilon)$ and hence the thesis follows.
Finally, direct computation shows that
$\Upsilon_\tau(\theta,\varepsilon,0) =
\Upsilon_\mu(\theta,\varepsilon,0)$ so that $
 \Delta^{\pm}_\tau(\theta,0) =   \Delta^{\pm}_\mu(\theta,0).$
\qed

\section{All Things Considered: Main Theorems, Comments and Examples}
\label{secallthings}

We can now complete our work, first by linking the inequalities
obtained so far and then by commenting on their optimality with
the aid of the von Neumann -- Kakutani Map \cite{von} and of an
intermittent map due to Pomeau -- Manneville \cite{pmv} and
Gaspard -- Wang \cite{gasp}. Recall that we have put ourselves in
a rather general setting, by requiring only ergodicity of the
dynamical system considered. Our fundamental result is therefore:
\begin{theorem}
When the dynamical system $(X,T,{\mathcal A},\mu)$ is ergodic and
$X$ is a compact metric space enclosed in ${\bf R}^n$, for any
$\theta \in \Theta$, the different dimensions defined in this work
are linked by the inequalities:
 \beq
 D^{\pm}_\tau(q) \leq \Delta^{\pm}_\tau(\theta,q)  \leq  \Delta^{\pm}_\mu(q)  =
 D^{\pm}_\mu(q), \;\;\mbox{ for $q>0$},
\nuq{ups9}
\beq
  \Delta^{\pm}_\tau(\theta,q)  \geq  \Delta^{\pm}_\mu(\theta,q)  \geq
 D^{\pm}_\mu(q), \;\;\mbox{ for $q \leq 0$},
\nuq{ups10} \beq
  \Delta^{\pm}_\tau(\theta,q)  \geq D^{\pm}_\tau(q), \;\;\mbox{ for $q \leq 0$},
\nuq{ups10b}
 and, for $q=0$,
\beq
  \Delta^{\pm}_\tau(\theta,0)  =
 \Delta^{\pm}_\mu(\theta,0).
 \nuq{ups10x}
 \label{teototu}
\end{theorem}

{\em Proof.} Use Lemma \ref{lemaboxc1} together with Propositions
\ref{prop2} and \ref{prop1}. \qed

The only equality that we have proven to hold in full generality
is between $\Delta^{\pm}_\tau(\theta,0)$ and
$\Delta^{\pm}_\mu(\theta,0)$, obviously when computed on the same
grid $\theta$. It is believed that $\Delta^{\pm}_\mu(\theta,0) =
D^\pm_\mu(0)$ should hold  in large generality \cite{tre}. When
this is the case, we can also assess that
$\Delta^{\pm}_\tau(\theta,0)$ does not depend on the grid $\theta$
and this provides us with a means of computing the capacity
dimensions $D^\pm_\mu(0)$ via return times. It is remarkable that
no exceptions to the desired equality were known until very
recently: that is, there exists a (still unpublished) case where
$\Delta^{-}_\mu(\theta,0)$ is strictly larger than $D^-_\mu(0)$
(S. Tcheremchantsev, private communication).

The situation occurring for $q>0$ is fully described by a single
chain of inequalities, (\ref{ups9}). We want now to show that they
can be strict. In fact, the return times dimensions
$D^{\pm}_\tau(q)$ may decay to zero when $q$ tends to infinity
even when measure dimensions do not. This can be regarded as a
consequence of ``short returns'', a rather general occurrence. In
fact, let
\begin{equation}
\rho(\varepsilon;m) := \mu (\{x \in X \mbox { s.t. }
\tau_{B_\varepsilon(x)} (x) = m \})
 \label{tdim5c}
 \end{equation}
be the distribution of the first return of a point $x$ into the
ball of radius $\varepsilon$ centered at $x$. Also consider the
integrated distribution $R(\varepsilon;k)$:
\begin{equation}
R(\varepsilon;k) := \sum_{m=1}^{k} \rho(\varepsilon;m).
 \label{tdim5d}
 \end{equation}
We have the following
\begin{lemma}
If for some $k \geq 1$, there exist constants $C$ and $\delta>0$
such that $R(\varepsilon;k) \geq C \varepsilon^\delta$, then
$D^{\pm}_\tau(q) \leq \frac{\delta}{q-1}$ for all $q>1$.
 \label{lemdeca1}
\end{lemma}
{\em Proof.} Let $q>1$. Clearly,
\begin{equation}
 \Gamma_\tau(\varepsilon,q) = \sum_{m=1}^{\infty} \rho(\varepsilon;m) m^{1-q}
  \geq k^{1-q} R(\varepsilon;k) \geq C k^{1-q} \varepsilon^\delta ,
 \label{nte1bx}
 \end{equation}
which yields the thesis. \qed

A similar Lemma holds obviously also for
$\Delta^{\pm}_\tau(\theta,q)$. This lemma shows that, roughly
speaking, in order for dimensions not to tend to zero when $q$
tends to infinity, the probability of small returns must vanish
faster than any power of $\veps$, when $\veps$ tends to zero. We
shall momentarily describe a system, the von Neumann -- Kakutani
map, where to the contrary this probability decays as $\veps$ and
the inequalities (\ref{ups9}) are strict for $q>2$.

In a previous work \cite{prl} we have outlined another mechanism
for short returns: the existence of fixed points of a continuous
map $T$. Indeed, $R(\varepsilon;1)$ can be bounded from below by
the measure of a box of radius proportional to $\varepsilon$
centered at any fixed point. This latter scales, for small
$\varepsilon$, with the local dimension at the fixed point, a
value that can be used in Lemma \ref{lemdeca1}. We must remark
that in the examples presented in \cite{prl} local dimensions
provide an upper bound, while the exact asymptotic result should
involve the correlation dimension $D_\mu(2)$ (see below).

Let us now consider the case $q < 0$. In full generality, we can
only establish the shorter chains of inequalities
(\ref{ups10},\ref{ups10b}). We are not able to perform other
comparisons. Contrary to what might seem at first blush, this is
not the result of a deficiency of our technique. In fact, as we
have remarked in Sect. \ref{sectriedi}, it may happen that
$\Delta^{+}_\mu(\theta,q)$ be larger than $D^{+}_\mu(q)$, even
infinite. This might happen not because of any peculiarity of the
measure, but because of the choice of the grid $\theta$. In turns,
this also ``spoils'' $\Delta^{+}_\tau(\theta,q)$, because of the
inequality (\ref{ups10}), but {\em not} $D^{+}_\tau(q)$, which is
{\em smaller than} $\Delta^{+}_\tau(\theta,q)$ and, as such, is
not linked to $\Delta^{+}_\mu(\theta,q)$.

\begin{figure}[ht]
\includegraphics[width=8cm,height=12cm,angle=270]{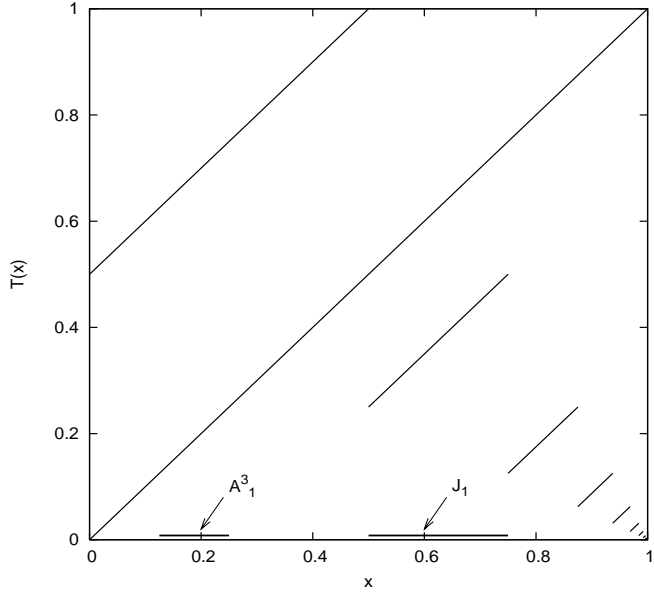}
\caption{Graph of the von Neumann -- Kakutani Map defined in eq.
(\ref{fze1}). Also drawn are the slope--one line through the
origin and the sets $J_1$ and $A^3_1$.} \label{figkak}
\end{figure}

This precisely happens for the von Neumann--Kakutani map
\cite{von}, described in detail in Sect. \ref{secmapintex} and
pictured in Figure \ref{figkak}, whose absolutely continuous
invariant measure is the uniform Lebesque measure over the unit
interval.
This map $T : [0,1] \rightarrow [0,1]$ is a sort of infinite
intervals exchange map that permutes diadic sub--intervals of any
order. In this permutation points in any binary interval of length
$2^{-n}$ (for any integer value of $n$) ``visit once'' all
remaining intervals before returning home: see Lemma
\ref{lemperm1} in Sect. \ref{secmapintex} below. From the point of
view of return times, this is a sort of dream situation, where all
points return in a time $\tau=2^n$ that is exactly the inverse of
the measure of the interval. As a consequence, for these sets,
formulae (\ref{hold2}) hold as equalities for all real values of
$s$. Nonetheless, Grassberger and Jensen {\em et al.} conjecture
is verified on partially for this dynamical system, as the
following theorem shows:

\begin{theorem}
In the dynamical system $(T,[0,1],\lambda)$, where $T$ is the von
Neumann -- Kakutani map defined in eq. (\ref{fze1}) and $\lambda$
is the Lebesgue measure, the dimensions defined in this work take
the values:
 \begin{equation}
  D^{\pm}_\tau(q) = \left\{
  \begin{array}{ll}
       1 & \mbox{ for } q \leq 2 \\
       \frac{1}{q-1} & \mbox{ for }  q > 2. \\
  \end{array} \right.
 \label{tdim4b}
 \end{equation}
Moreover, for all grids $\theta$,
 \beq
\Delta_\tau^+(\theta,q) = \Delta_\mu^+(\theta,q) = \infty
 \;\; \mbox{ for } q < 0.
 \nuq{td4b3}
Finally, there are an infinite number of grids $\theta$ for which
\begin{equation}
  \begin{array}{ll}
      \Delta_\tau^+(\theta,q) = 1 & \mbox{ for } q \geq 0,  \\
     \Delta_\tau^-(\theta,q) = 1 & \mbox{ for } q \leq 0.    \\
  \end{array} 
 \label{tdim4b2}
 \end{equation}
\label{teokaku}
\end{theorem}
 {\em Proof.} See Sect. \ref{secprteokaku}.

Recall now that for the Lebesgue measure on the unit interval, for
any $\theta$, we have that $\Delta^{-}_\mu(\theta,q) =
D^{\pm}_\mu(q) = 1$ for any $q$, $\Delta^{+}_\mu(\theta,q) = 1$
for $q \geq 0$ and $\Delta^{+}_\mu(\theta,q) = \infty$ for $q <
0$. Therefore, in this case, measure and return time generalized
dimensions $D_\mu^\pm(q)$ and $D_\tau^\pm(q)$ coincide for all $q
\leq 2$. At the same time, for $q<0$, $\Delta_\tau^+(\theta,q)$
and $\Delta_\mu^+(\theta,q)$ are affected by the ``edge effect''
discussed above and feature a ``pathological'' value.

As a consequence of the short--returns phenomenon discussed
earlier in this section, Lemma \ref{lemdeca1}, $D_\mu^\pm(q)$ and
$D_\tau^\pm(q)$ differ for $q > 2$, and the latter dimensions
vanish for large $q$. Observe also a ``phase transition'' behavior
occurring at $q=2$. Remark finally that, by choosing particular
grids, we obtain equality also for the grid dimensions, when
taking the superior limit (for $q>0$) and the inferior limit in
the opposite case.

We conjecture that what observed for this map is a rather common
situation: that is to say, we expect that
\begin{conjecture}
\label{conje}  For a large class of dynamical systems
$D_\tau^\pm(q)= D_\mu^\pm(q)$ for $q_c < q \leq 2$ ($q_c$ being
the lowest value of $q$ for which partition functions of return
times are finite, recall Lemma \ref{lem-existence} and see below
for an example) and $D_\tau^\pm(q)= D^\pm_\mu(2)/(q-1)$ for $ q
\geq 2$ (exactly, or at least asymptotically for large $q$).
\end{conjecture}

At this point, it is relevant to quote the results of
\cite{ugacha} that have already been mentioned in Sect.
\ref{secintro}. They hold under strong assumptions on the
dynamical system under investigation and for cylinders rather than
balls ({\em i.e.} without relation to the distance function). In
fact, it has been shown that for Bowen--Gibbs measures, defining
partition functions and generalized dimensions for entrance
(rather than return) times in dynamical cylinders, these latter
coincide with Renyi entropies for $q<2$ and behave as
$P(2\phi)/(q-1)$ for larger $q$. Here $P$ is the topological
pressure of the potential $\phi$ defining the Bowen--Gibbs
measure. For return times in cylinders, though, only the statement
for $q<1$ has been derived. These results outline interesting
techniques that might possibly be improved, and complemented with
geometric considerations, to prove in vast generality the
relations between measure and return time dimensions for balls, as
originally conjectured in \cite{prl} and formulated above in a
more precise form.

\begin{figure}[ht]
\includegraphics[width=8cm,height=12cm,angle=270]{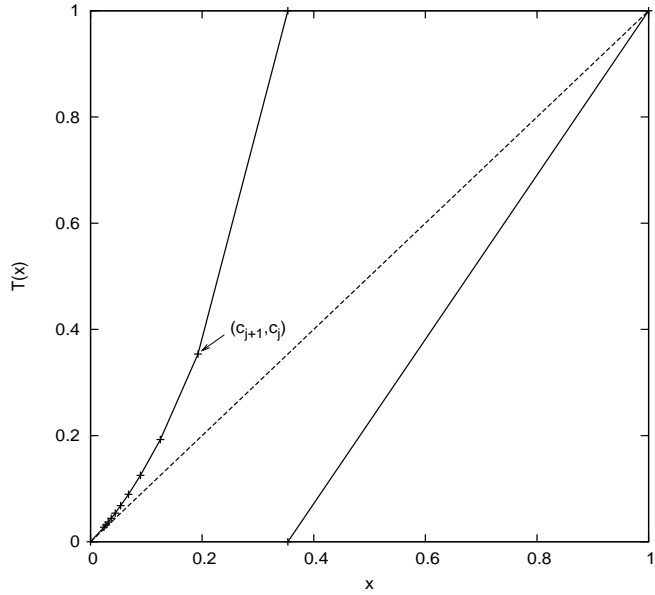}
\caption{Graph of the Gaspard-Wang Map defined in eq.
(\ref{mapxav2}), for $p=-3/2$, see Sect. \ref{secintermap} for
details. Also drawn is the slope--one line through the origin.}
\label{figwang}
\end{figure}

Let us now move to a final example, which shows that inequalities
between measure and return time dimensions may be strict also for
negative values of $q$: in fact, $D^{\pm}_\tau(q)$ and
$\Delta^{\pm}_\tau(q)$ may be infinite for all $q$ smaller than a
critical value  $q_c < 0$, when $D^{\pm}_\mu(q)$ is finite. This
is notably the case of intermittent maps, the simplest of which is
perhaps the Gaspard -- Wang \cite{gasp} piece--wise linear
approximation of the Pomeau Manneville map \cite{pmv} described in
Sect. \ref{secintermap} and pictured in Fig. \ref{figwang}. This
is a map of the unit interval into itself, with an absolutely
continuous invariant measure. Zero is a fixed point of the map and
the dynamics may spend arbitrarily long time spans in its
neighbourhood. For this dynamical system we can prove the
following theorem, that demonstrates a case where $D^{\pm}_\tau(q)
> D^{\pm}_\mu(q)$ for sufficiently negative $q$, an inequality
that is specific to this particular case and is not included among
those in formulae (\ref{ups10},\ref{ups10b}).

\begin{theorem}
In the dynamical system $(T,[0,1],\mu)$, where $T$ is the Gaspard
-- Wang map defined in eq. (\ref{mapxav2}) with parameter $p<-1$
and $\mu$ is its unique absolutely continuous invariant measure,
return time dimensions satisfy $D^\pm_\tau(q) =
\Delta^{\pm}_\tau(q) = \infty$ for all $q < q_c := p+1$, while
measure dimensions take the values $D^\pm_\mu(q) = 1$ for $q \leq
-p$ and $D^\pm_\mu(q) = (1+ \frac{1}{p}) \frac{q}{q-1}$ for $q
\geq -p$.
 \label{teomane}
\end{theorem}

{\em Proof.} See Sect. \ref{secintermap}.

\section{Conclusions}
\label{secconcl} We might now try to conclude by saying that the
idea to use return times in a straightforward way to compute
generalized measure dimensions, following the programme whose
history has been briefly outlined in the Introduction, is only
applicable after a detailed analysis of the dynamical system
considered.

The general inequalities that we have derived clarify the mutual
relations among the dimensions that we have defined. As a
by--product, these inequalities provide {\em universal bounds} for
the global statistics of return times that hold for {\em all}
ergodic dynamical systems possessing a given invariant measure
$\mu$, and indeed also for a large class of stochastic processes
having invariant distribution $\mu$.

We have found examples where $D_\mu(q)$ and $D_\tau(q)$ differ for
$q>2$, or for $q<q_c$. At the present moment, we do not know of
any example where $D_\mu(q)$ and $D_\tau(q)$ differ in the
interval $(q_c,2)$. Nevertheless, we are not able to prove
equality in full generality with the means employed in this paper.
We consider this, as well as the precise formulation and proof of
conjecture \ref{conje}, to be a point of utmost interest for
future investigations.

Turning from the general case to specific applications, we feel
that one could prove part or all of Conjecture \ref{conje} with
problem--specific tools. This might indeed be good news, that
would partly fulfill the original Grassberger  and Jensen {\em et
al.} program, in particular for dimensions with negative $q$, that
are known to be more elusive to compute numerically and more
intriguing theoretically than those for positive $q$
\cite{due,sator}.

Finally, whether linked to generalized measure dimensions or not,
the moments of return times studied in this work deserve attention
in their own, in our view. In fact, at difference with local
quantities studied in the literature (such as probabilities of
return to shrinking neighborhoods of a given point---a well
examined topic, see {e.g.} \cite{ornstein,beno,beno2,nico,miguel})
they provide a global characteristic of the dynamics of a system.

The remainder of this paper consists now of three sections giving
details and proofs for the two maps quoted in this paper.

\section{The Map of von Neumann and Kakutani}
\label{secmapintex}

In this section we present the details of the intervals exchange
map due to von Neumann and Kakutani \cite{von}, mentioned in Sect.
\ref{secallthings}. The basic properties of this map are known,
but we prefer to re-derive them here for completeness and because
they help us to understand some subtleties of return times  for
this map.

We start by defining two families of intervals in $[0,1]$. The
first is
\begin{equation}
J_k := [1-2^{-k},1-2^{-k-1}),  \; k=0,1,\ldots \label{int1}
 \end{equation}
Clearly, $X = [0,1] = \bigcup_{n=0}^\infty I_n \bigcup \; \{1\}$.
Then, the map of von Neumann and Kakutani, $T$, is defined as
follows:
\begin{equation}
 T(x) := \left\{
  \begin{array}{ll}
       x - 1  + 2^{-k} + 2^{-k-1} & \mbox{ for } x \in J_k \\
      0 & \mbox{ for } x = 1 \\
  \end{array} \right.
 \label{fze1}
 \end{equation}
The map $T$ is piece-wise continuous, composed of an infinite
number of affine segments and invertible (except for the point
$x=1$ which has no preimage).

In addition,  for any positive integer $n$, define a measurable
partition of $X$ in (open) binary intervals:
\begin{equation}
A^n_j := (j 2^{-n}, (j+1) 2^{-n}),  \; j=0,\ldots,2^n-1.
 \label{int2}
 \end{equation}
All but a finite number of points in $X$ are covered by the
partition. Exception are the {\em boundary points} $\zeta_k^n = k
\; 2^{-n}$, with $k = 0,\ldots,2^n$.
 Then, it is easy to see that
\begin{lemma}
 For any positive $n$, the map $T$ permutes the family of
intervals $\{A^n_j\}$. The permutation is cyclic, of period
$N=2^n$, in the in the sense that $T^N(A^n_j) = A^n_j$ for any
$j$, and no shorter $N$ exists with this property.
\label{lemperm1}
\end{lemma}
{\em Proof.} Let $\sigma = \sigma_1 \sigma_2 \ldots \sigma_n$ be
the binary expansion of $j$, defined as follows (notice the order
of digits):
\begin{equation}
   j := \sum_{k=1}^{n} \sigma^k 2^{n-k}.
  \label{fj}
 \end{equation}
Let also $k(\sigma)$ be the index of the first zero in $\sigma$:
 \begin{equation}
 k(\sigma) := \min \;  \{ j \in {\bf N} \mbox{ s.t. } \sigma_j = 0 \}.
 \label{fzero}
 \end{equation}
Intervals $A^n_j$ shall therefore be labelled as $A^n_\sigma$,
where $\sigma$ is a word of length $n$. We shall also use the
complementary digit function $\bar{\cdot}$, where $\bar{0}=1$,
$\bar{1}=0$.

All points in the interval $A^n_\sigma$ can be written in binary
form as $x = 0. \omega_1 \ldots \omega_n \omega_{n+1} \ldots$,
where $\omega_i = \sigma_i$ for $i=1,\ldots,n$ and where
$\omega_{n+1}, \ldots$ is any infinite sequence of digits (except
for the sequence composed of all ones). It can be verified that,
in binary notation, the map $T$, eq. (\ref{fze1}), corresponds to
the symbolic map $S(\omega) = . \eta_1 \eta_2 \ldots$, with
 \begin{equation}
 \eta_j := \left\{
  \begin{array}{ll}
       \overline{\omega_j} & \mbox{ for } j \leq k(\omega) \\
       \omega_j & \mbox{ for } j > k(\omega) \\
  \end{array} \right.
 \label{fzero2}
 \end{equation}
Therefore, any interval $A^n_\sigma$ is mapped into the interval
$A^n_\eta$, labelled by the first $n$ digits of $\eta$. For this
reason, with a slight misusage of notation, we shall indicate by
$S$ also the map $\sigma \rightarrow \eta$ on the set of
$n$-letter words, or equivalently via eq. (\ref{fj}) on the set of
integers $[0,2^n-1]$. The map $S$ acts a cyclic permutation of all
intervals $A^n_\sigma$, of period $N=2^n$, for any value of $n$,
the length of the word $\sigma$. \qed

An interesting consequence of the previous lemma is the following
proposition,
\begin{proposition}
The Lebesgue measure $\lambda$ on $X$ is invariant and ergodic for
the action of the map $T$. Moreover, the dynamical system
$(X,T,\lambda)$ is metrically and topologically transitive, but
not mixing. \label{proptmap}
\end{proposition}
{\em Proof.} The first statement is almost immediate from the form
of the map $T$,  eq. (\ref{fze1}) and the first part of Lemma
\ref{lemperm1}: given any open interval $I$, its counter-image is
a finite union of disjoint intervals whose lengths add up to the
length of $I$.

To prove ergodicity one needs to show that for any measurable sets
$B$ and $C$,
\begin{equation}
  \lim_{m \rightarrow \infty} \frac{1}{m} \sum_{k=0}^{m-1}
   \mu( T^{-k}(B) \cap C  ) = \mu(B) \mu(C).
  \label{fcard0}
 \end{equation}
Let $B$ and $C$ be finite unions of binary intervals $A^n_j$ at
resolution $n$. Indicate with $\mathcal B $ and $\mathcal C $ the
sets of indices of the intervals composing the sets $B$ and $C$,
like in $
  B := \bigcup_{j \in \mathcal B} A^n_j.
 $
Finally let $\#(\mathcal B)$ and $\#(\mathcal C)$ be the
cardinalities of these sets, respectively. Recall that $T$
permutes the intervals $A^n_j$ as in $A^n_j$ as in Lemma
\ref{lemperm1} and so does $T^{-1}$. Therefore,
\begin{equation}
   T^{-k}(B) =  T^{-k}(\bigcup_{j \in \mathcal B} A^n_j) =
   \bigcup_{j \in \mathcal B} T^{-k}(A^n_j),
  \label{fcard3}
\end{equation}
and the intervals in the union above are disjoint, so that
\begin{equation}
  \mu( T^{-k}(B) \cap C ) =
   \sum_{j \in \mathcal B} \mu (T^{-k}(A^n_j) \cap C).
  \label{fcard3b}
\end{equation}

Moreover, for any $j$ and $k$, either $T^{-k}(A^n_j)$ has empty
intersection with $C$, or it coincides with one of the binary intervals
composing $C$.
 It is then convenient to just consider the interval index
map, that we have also indicated by $T$. Define therefore the set of
``times'' for which such intersection is not empty:
\begin{equation}
 N^n_{\mathcal C}(j) :=
 \{ k \in {\bf Z} \mbox{ s.t. } 0 \leq k < 2^n \mbox{ and } T^{-k}(j) \in {\mathcal C} \}.
 \label{fcard1}
 \end{equation}
For each $k \in  N^n_{\mathcal C}(j)$ we have that
 \begin{equation}
\mu( T^{-k}(A^n_j) \cap C ) = \mu( T^{-k}(A^n_j)) = \mu(A^n_j) = 2^{-n},
 \label{fcard2b}
 \end{equation}
and for $k \not \in  N^n_{\mathcal C}(j)$, $\mu( T^{-k}(A^n_j)
\cap C ) =0$. We then compute
\begin{eqnarray}
  \sum_{k=0}^{2^n-1} \mu( T^{-k}(B) \cap C ) =
   \sum_{j \in \mathcal B} \sum_{k=0}^{2^n-1} \mu (T^{-k}(A^n_j) \cap C) =
 \nonumber \\ =
\sum_{j \in \mathcal B} \sum_{k \in N^n_{\mathcal C}(j)}
\mu(A^n_j) = 2^{-n} \sum_{j \in \mathcal B} \# N^n_{\mathcal
C}(j).
  \label{fcard4}
\end{eqnarray}

Finally, observe that being $T$ a cyclic permutation of the first
$2^n$ integers, $T^{-k}(j) \in \mathcal C$ holds $\#(\mathcal C)$
times along any cycle of times of length $2^n$:
 \begin{equation}
   \# N^n_{\mathcal C}(j) = \#(\mathcal C),
 \label{fcard1b}
 \end{equation}
which means
 \begin{equation}
 \frac{1}{2^n}
\sum_{k=0}^{2^n-1} \mu( T^{-k}(B) \cap C ) = 2^{-2n} \#(\mathcal
C) \#(\mathcal B) = \mu(B) \mu(C).
 \label{fcard5}
 \end{equation}

This easily entails the limit (\ref{fcard0}) for binary intervals.
Since these latter generate the Borel sigma algebra, the result
follows generally.

Topological transitivity (ergodicity) is easily implied by Lemma
\ref{lemperm1}, since given any two open sets $B$ and $C$ there
exist  an $n > 0$ and $A^n_j$, $A^n_{j'}$, $0 \leq j,j' < 2^n$,
such that $A^n_j \subset B$ and $A^n_{j'} \subset C$. Choose then
$k$ such that $T^k(j)=j'$ to obtain the result.

It is also immediate to see that strong mixing is not present: in
fact, this is ruled out by the cyclic nature of the images
$T^{-k}(A^n_j))$, for a single $n,j$.

Weak mixing can be ruled out by a careful usage of eq.
(\ref{fcard5}). \qed

The interesting properties of the map $T$ so defined permit us to
prove the following
\begin{proposition}
In the dynamical system $(T,\lambda)$ defined in this section,
over the sequence $\varepsilon_n = 2^{-n}$, the return time
partition function $\Psi_\tau(\theta,\varepsilon,q)$ can be
explicitly computed when $\theta$ is the grid having origin at
zero. \label{proprho}
\end{proposition}

{\em Proof.} Let $\varepsilon_n = 2^{-n}$ and consider the newly
introduced partition function $\Psi_\tau(\theta,\varepsilon,q)$,
eq. (\ref{nted3}).
It requires the computation of $\tau_{\overline{A^n_j}} (x)$ for
$x \in A^n_j$. This we shall do now.

Observe first that $\overline{A^n_j} = \bigcup_{i=-1,0,1}
A^n_{j+i}$, where obviously $A^n_j := \emptyset$ for $j<0$ or $j
\geq 2^n$. Because of the Lemma \ref{lemperm1}, in the von Neumann
-- Kakutani map $T$, these times are independent of the point $x$
in $A^n_j$ and can be computed in terms only of the index map $S$.
Let $j=j(x)$ the index of the binary interval containing $x$.
Then,
\begin{equation}
  \tau_{\overline{A^n_j}} (x) = \inf \{ k \geq 1 \mbox{ s.t. }
   S^k(j) \in \{j-1,j,j+1\} \}.
 \label{tdim3}
 \end{equation}

It is now possible (although complicated) to compute {\em
explicitly} the distribution of first return times of the index
$j$ into $\{j-1,j,j+1\}$. An example will be provided in the next
section. Observe that these return times over a finite set can
take any value between one and the cardinality of the finite set,
$2^n$. Let $\varrho^n_m$ the cardinality of the set of values of
$j$ for which the return value is $m$. Then, the following formula
holds:

\begin{equation}
  \varrho^n_m = \left\{
  \begin{array}{ll}
       2^l & \mbox{ for } m = 3 \times 2^l, \; l=0,\ldots,n-3 \\
       2^{n-2} & \mbox{ for }  m = 2^{n-2} \\
        2^{n-1} + 1 & \mbox{ for }  m = 2^{n-1} \\
                0 & \mbox{ elsewhere} \\
  \end{array} \right.
 \label{tdim4}
 \end{equation}

Therefore, letting
\begin{equation}
\rho^{(n)}(m) := \mu (\{x \in X \mbox { s.t. } \tau_{\overline{A^n_j}} (x)
= m \}) = 2^{-n} \varrho^n_m
 \label{tdim5}
 \end{equation}
we can define a family of distribution functions over the integers larger
than, or equal to one, according to which

\begin{equation}
 \Psi_\tau(\theta,\varepsilon,q) := \sum_j
  \int_{A_j} \tau_{\overline{A}_j}^{1-q}(x) \; d\mu(x) =
  \sum_{m=1}^\infty m^{1-q} \rho^{(n)}(m),
 \label{nted3b}
 \end{equation}
with $\varepsilon = 2^{-n}$. \qed

\section{Proof of Theorem \ref{teokaku} on The Map of von Neumann and Kakutani}
\label{secprteokaku}

Consider the measurable partitions $\{A^n_j\}$ of $X$ in binary
intervals of length $2^{-n}$ defined in eq. (\ref{int2}) in Sect.
\ref{secmapintex} and let $A^n_x$ be the element of the partition
containing the point $x$. We call $n$ the order of the partition.
Lemma \ref{lemperm1} implies that
 \begin{equation}
\tau_{A^n_x}(x) = 2^n,
 \label{kaku1}
\end{equation}
for any integer $n$ and for any $x \in X$. We part the proof of
Theorem \ref{teokaku} in several sections.

\subsection{Part a, where large balls completely cover dyadic intervals.}

In fact, when $2^{-n} \leq \veps \leq 2^{-n+1}$ the ball of radius
$\veps$ centered at $x$ covers the dyadic interval including $x$:
$B_\varepsilon(x) \supset A^n_x$. Therefore,
 \begin{equation}
\tau_{B_{ \varepsilon}(x)} (x) \leq \tau_{A^n_x}(x) = 2^n.
 \label{kaku2}
\end{equation}
 Choose $q$ such that $1-q \; \gtle \; 0$. We now derive
inequalities bounding the partition function
$\Gamma_\tau(\varepsilon,q)$ in both cases. Firstly,
\begin{equation}
\int  [\tau_{B_\varepsilon(x)}(x)]^{1-q}
 d\mu(x) \legt
\int  [\tau_{A^n_x}(x)]^{1-q}
 d\mu(x) =
 2^{n(1-q)}.
 \label{kaku3}
\end{equation}
Observe that substituting the inequalities linking $\veps$ and $n$
one gets
\begin{equation}
 \log \Gamma_\tau(\varepsilon,q) \leq (q-1) \log \varepsilon
 \label{kaku4}
\end{equation}
for $1-q > 0$ and
\begin{equation}
 \log \Gamma_\tau(\varepsilon,q) \geq (q-1) (\log \varepsilon +
 \log 2)
 \label{kaku5}
\end{equation}
for $1-q < 0$. This implies that
\begin{equation}
 D^\pm_\tau(q) \leq 1
 \label{kaku6}
\end{equation}
for all values of $q \neq 1$ (A similar treatment could also yield
the case $q=1$, we do not include it here for conciseness).

\subsection{Part b, where we exploit balls included in dyadic intervals of
order $n$.} Let $\veps$ and $n$ be related as:
\begin{equation}
\frac{1}{8} 2^{-n} \leq \veps \leq \frac{1}{4} 2^{-n}.
 \label{kaku6b}
\end{equation}
Then, when $x$ is close to the midpoint of the interval $A^n_x$,
$B_\varepsilon(x) \subset A^n_x$, so that
 \begin{equation}
\tau_{B_{ \varepsilon}(x)} (x) \geq \tau_{A^n_x}(x) = 2^n.
 \label{kaku7}
\end{equation}
The same inequality also holds for points close to zero and one,
the extrema of $X$, because in this case $(X \cap
B_\varepsilon(x)) \subset A^n_x$. Collectively, all these points
define the set $G_{n,\veps}$. It is easy to see that because of
eq. (\ref{kaku6b}) the measure of this set amounts to at least
half of the total measure. We now take $1-q > 0$, so that
\begin{equation}
\int_X  [\tau_{B_\varepsilon(x)}(x)]^{1-q}
 d\mu(x) \geq \int_{G_{n,\veps}}  [\tau_{B_\varepsilon(x)}(x)]^{1-q}
 d\mu(x) \geq
\int_{G_{n,\veps}}  [\tau_{A^n_x}(x)]^{1-q}
 d\mu(x) \geq \frac{1}{2}
 2^{n(1-q)}.
 \label{kaku8}
\end{equation}
Proceeding as above, we find that
\begin{equation}
 D^\pm_\tau(q) \geq 1
 \label{kaku9}
\end{equation}
for all values of $q \leq 1$. Together with eq. (\ref{kaku6}) this
implies that
\begin{equation}
 D^\pm_\tau(q) = 1
 \label{kaku10}
\end{equation}
for all values of $q < 1$.

\subsection{Part c, where we exploit balls included in dyadic intervals of orders $0$
to $n$.}
 Let the inequalities (\ref{kaku6b}) still hold. We extend
the argument of part b. Suppose that $x$ does not belong to
$G_{n,\veps}$. This means that $x$ is within $\veps$ of any of the
endpoints of the interval $A^n_x$ internal to $[0,1]$. Then, the
ball $B_\varepsilon(x)$ is not included in $A^n_x$, but it
stretches to reach one neighboring element of the partition.

Forcefully, $B_\varepsilon(x)$ includes a boundary point of the
measurable partition of order $n$, of the form $\zeta_k^n = k \;
2^{-n}$, with $k$ integer. Clearly, two cases are possible: either
$\zeta_k^n$ is an ``even'' point ($k$ even, in which case it is
also a boundary point of the partition of order $n-1$), or it is
an ``odd'' point. If it is an odd point, $B_\varepsilon(x)$ is
necessarily included in $A^{n-1}_x$ and $\tau_{B_{
\varepsilon}(x)} (x) \geq \tau_{A^{n-1}_x}(x) = 2^{n-1}$.  Let
$G_{\veps,n-1}$ be the set of points $x$ that are $\veps$--close
to odd boundary points of the partition of order $n$. It is
immediate that this set is composed of $2^{n-1}$ intervals of
length $2 \veps$.

The construction can be clearly iterated, by considering among
even boundary points of level $n$ those which are odd at level
$n-1$: this defines a new set of points $G_{\veps,n-2}$ consisting
of $2^{n-2}$ intervals of length $2 \veps$. For any $x$ belonging
to this set, $\tau_{B_{ \varepsilon}(x)} (x) \geq
\tau_{A^{n-2}_x}(x) = 2^{n-2}$ and so on. The last stage of the
construction is the set $G_{\veps,0} =
(\frac{1}{2}-\veps,\frac{1}{2}+\veps)$, for which, rather
trivially, $\tau_{B_{ \varepsilon}(x)} (x) \geq \tau_{A^{0}_x}(x)
= 2^{0}=1$.

The above proves the following formul\ae: for any $\veps > 0$ and
$n$ satisfying the inequalities (\ref{kaku6b}), one has:
\begin{equation}
 X = \bigcup_{j=0}^n G_{\veps,j},
 \label{kaku13}
\end{equation}
with $\lambda(G_{\veps,j} \bigcap G_{\veps,j'}) = 0$ if $j \neq
j'$ (recall that $\lambda$ is the Lebesgue measure and that
$G_{\veps,n}$ has been defined in point b), so to provide another
measurable partition of $X$. Moreover,
\begin{equation}
 2^{-n+j} \geq  2^{-n-1+j} \geq  \lambda(G_{\veps,j}) \geq
2^{-n-2+j},
 \label{kaku14}
\end{equation}
for $0 \leq j \leq n-1$ and $1 \geq \frac{3}{4} + 2^{-n-2} \geq
\lambda(G_{\veps,n}) \geq 2^{-1}$. Finally,
\begin{equation}
 \tau_{B_{
\varepsilon}(x)} (x) \geq 2^{j}
 \label{kaku15}
\end{equation}
for any $x \in G_{\veps,j}$ and $j = 0,\ldots,n$.

We can now evaluate the partition function: let $1-q < 0$, {\em
i.e.} $q>1$, so that
\begin{equation}
\Gamma_\tau(\varepsilon,q) = \sum_{j=0}^n \int_{G_{\veps,j}}
[\tau_{B_\varepsilon(x)}(x)]^{1-q}
 d\mu(x) \leq
 \sum_{j=0}^n \lambda(G_{\veps,j}) 2^{j(1-q)}
\leq  \sum_{j=0}^n 2^{-n+j} 2^{j(1-q)},
 \label{kaku16}
\end{equation}
where we have used the widest inequality in (\ref{kaku14}) to
obtain a simpler formula. In fact, eq. (\ref{kaku16}) easily
yields
\begin{equation}
\Gamma_\tau(\varepsilon,q)  \leq G(q,n) := 2^{-n}
\frac{2^{(2-q)(n+1)} -1}{2^{(2-q)} - 1},
 \label{kaku17}
\end{equation}
where the function $G(q,n)$ has been defined. Two cases must now
be considered, in the asymptotics of $G(q,n)$ as $n$ tends to
infinity, or $\veps$ goes to zero, according to (\ref{kaku6b}).
First, when $2>q>1$ we find $log(G(q,n)) \sim \veps^{q-1}$, so
that $D^\pm_\tau(q) \geq 1$ and finally
\begin{equation}
 D^\pm_\tau(q) = 1
 \label{kaku18}
\end{equation}
for all values of $q \leq 2$. in the other case, $q>2$, we find
$log(G(q,n)) \sim \veps^{-1}$, so that
\begin{equation}
 D^\pm_\tau(q) \geq \frac{1}{q-1}.
 \label{kaku19}
\end{equation}

\subsection{Part d, where we exploit the return
properties of a particular sequence of dyadic intervals.}
This is the last part of this proof.  In Proposition \ref{proprho}
we have computed the return times of dyadic intervals in their
enlarged neighborhood, according to sections \ref{sectriedi},
\ref{bvcrtdim}. We need a particular case of Prop. \ref{proprho}
that can be easily proven. It appears from eq. (\ref{tdim4}) that
{\em for any} $n$ there exist {\em one} interval $A^n_j$ with
return time 3 inside $\overline{A}^n_j$. The index of this
interval is $j=2^{n-1}-1$ and its symbolic address is $ 0 1 \ldots
1$. Let this value of $j$ be fixed in the following. It maps to $1
1 \ldots 1$ and successively $0 \ldots 0$ and $1 0 \ldots 0$. It
is easy to see that the first and the last are neighboring
intervals. In coordinates, the quasi-cycle is $A^n_j =
(\frac{1}{2} - 2^{-n},\frac{1}{2}) \rightarrow (1 - 2^{-n},1)
\rightarrow (0 ,2^{-n}) \rightarrow (\frac{1}{2},\frac{1}{2} +
2^{-n}) = A^n_{j+1}$. It contains the orbit  $1 \rightarrow 0
\rightarrow \frac{1}{2}$. The sequence $1 \rightarrow 0
\rightarrow \frac{1}{2} \rightarrow 1 $  does {\em not} belong to
any orbit, but it is arbitrarily well approximated by true orbits.

If we now take  $2^{-n+2} \geq \veps \geq 2^{-n+1}$ we have that
$B_\varepsilon(x) \supset A^n_{j+1}$ for any $x$ in $A^n_j$ and
the above implies that $\tau_{B_\varepsilon(x)}(x) \leq 3$ on
$A^n_j$. Therefore, when $1-q<0$,
\begin{equation}
\int_X  [\tau_{B_\varepsilon(x)}(x)]^{1-q}
 d\mu(x) \geq
 \int_{A^n_j}  [\tau_{B_\varepsilon(x)}(x)]^{1-q}
 d\mu(x) \geq 3^{1-q} 2^{-n}.
 \label{kaku20}
\end{equation}
This yields $D^\pm_\tau(q) \leq \frac{1}{q-1}$ for any $q>1$ and
finally
\begin{equation}
 D^\pm_\tau(q) = \frac{1}{q-1}
 \label{kaku21}
\end{equation}
for any $q \geq 2$. This ends the proof of the part of the theorem
concerning $D^\pm_\tau(q)$.

\subsection{Proof of the results for box dimensions.}
As for the box dimension $\Delta^{\pm}(\theta,q)$, choose now the
grid $\theta$ with origin at zero and consider the sequence
$\veps_n = 2^{-n}$. Recall equation (\ref{kaku1}). It implies that
$ \Upsilon_\tau(\theta,\veps_n,q) = 2^{n(1-q)} = \veps_n^{q-1}$.
Then, for all $q \neq 1$
\begin{equation}
     \lim_{n \rightarrow \infty}
  \frac{1}{q-1}
   \frac{\log \Upsilon_\tau(\theta,\veps_n,q)}
  {\log \veps_n} = 1.
 \label{garfi0}
 \end{equation}
 Let now $q<0$. Clearly,
\begin{equation}
     \lim_{n \rightarrow \infty}
  \frac{1}{q-1}
   \frac{\log \Upsilon_\tau(\theta,\veps_n,q)}
  {\log \veps_n} \geq
  \liminf_{\veps \rightarrow 0}
  \frac{1}{q-1}
   \frac{\log \Upsilon_\tau(\theta,\veps,q)}{\log \veps} :=
  \Delta^{-}_\tau(\theta,q) \geq \Delta^{-}_\mu(\theta,q) = 1,
 \label{garfi1}
 \end{equation}
where the second inequality is eq. (\ref{ups10b}) and where the
last equality can be easily obtained. Since the first limit exists
and is equal to one, eq. (\ref{garfi0}), so is
$\Delta^{-}_\tau(\theta,q)$. For negative $q$ a similar argument
applies, which now requires the superior limit. The same results
are obviously found when the origin of the grid $\theta$ is a
point of the form $k \, 2^{-m}$, with integer $k$ and $m$. This
ends the proof of the Theorem.
 \qed

\section{The map of Gaspard and Wang}
\label{secintermap}

We now describe the piece--wise linear approximation of the Pomeau
Manneville intermittent map \cite{pmv} due to Gaspard and Wang
\cite{gasp}. Let $\{c_j\}_{j \in {\bf N}}$ be an ordered,
decreasing sequence of numbers between zero and one such that
$c_0= 1$ and such that $c_j$ tends to zero as $j$ tends to
infinity. Let $I_j := (c_{j+1},c_j)$ be the elements of a
partition of $[0,1]$ into open intervals of length $l_j =
c_j-c_{j+1}$. The map $T$ is defined as the transformation which
maps affinely and with positive slope $I_j$ onto $I_{j-1}$ for $j
\geq 1$ and $I_0$ onto $[0,1]$:
\begin{equation}
 T(x) = (x-c_{j+1})\frac{l_{j-1}}{l_j}  + d_j,
 \label{mapxav2}
 \end{equation}
for $x \in I_j$ and where we set $l_{-1}=1$ and $d_j=c_j$ for $j >
0$, $d_0=0$. Because this behavior, it is also called an {\em
infinite renewal chain} map. We now chose a family of such maps,
parameterized by $p < -1$, for which
 \begin{equation}
 c_j = (j+1)^p.
 \label{xav3}
 \end{equation}

{\em Proof of Theorem \ref{teomane}}. The absolutely continuous
invariant measure on $[0,1]$, whose density is constant on each
$I_j$ can be easily computed. One finds
\begin{equation}
\mu(I_j) = c_j a
 \label{xav2}
 \end{equation}
where the parameter $a=\mu(I_0)$ can be chosen so to normalize the
measure, of course when the sequence $\{c_j\}_{j \in {\bf N}}$ is
summable, which is always the case when $p<1$. The motion of this
dynamical system is such that $I_0$ can be also parted into an
infinity of adjacent intervals $K_j$, $j=0,\ldots,$ whose points
return to $I_0$ after exactly $j+1$ steps: if $x \in K_j$,
$\tau_{I_0}(x)=j+1$. One finds easily that
\begin{equation}
T(K_j)=I_j,
 \label{xav4}
 \end{equation}
so that all $K_j$ can be obtained by an affine transformation of
$I_j$: $K_j = (1-c_1) I_j + 1$. The measure of $K_j$ is
proportional to its length and hence to the length of $I_j$:
\begin{equation}
    \mu(K_j) = a l_j.
 \label{xav5}
 \end{equation}
Therefore, not all moments of the return times of points of $I_0$
into itself are finite: in fact,
\begin{equation}
  \int_{I_0} \tau_{{I}_0}^{1-q}(x) \; d\mu(x) = \sum_{j=0}^\infty
  \int_{K_j} \tau_{{I}_0}^{1-q}(x) \; d\mu(x) =
  a  \sum_{j=0}^\infty (j+1)^{1-q} l_{j}.
 \label{xav6}
 \end{equation}
Since $l_j \sim j^{p-1}$, the above integral is convergent only
when $q > p+1$. Also observe that formulae (\ref{hold2}) are here
strict inequalities.

Now, let us cover the unit interval by a box grid of side $\veps$
and let's evaluate the partition function
$\Psi_\tau(\theta,\varepsilon,q)$. Let us consider the particular
box $A$ that contains the point $c_1$ in its interior. If $c_1$ is
a boundary point, consider the box whose left extremum is $c_1$.
For any $\veps>0$ the box $A$ contains an infinite number of
$K_j$, those with $j>j_\veps$. At the same time, the enlarged box
$\overline{A}$ is enclosed in the union $I_0 \cup I_1 \cup \ldots
I_m$, when $m$ depends on $\veps$. If $\veps$ is sufficiently
small, we can take $m=1$. Since $K_j$ maps to $I_j$ and this to
$I_{j-1}$, et cetera, the time to return to $\overline{A}$ is
larger that the time to enter $I_m$ from $K_j$. This latter is,
clearly, $j-m$. Entering this information in eq. (\ref{nted3}) we
find that, when $q<1$,
\begin{equation}
 \Psi_\tau(\theta,\varepsilon,q) \geq
  \int_{A} \tau_{\overline{A}}^{1-q}(x) \; d\mu(x)
  \geq \sum_{j > j_\veps} \mu(K_j) (j-m)^{1-q}.
 \label{xav7}
 \end{equation}
Therefore, as in eq. (\ref{xav6}), the sum diverges for $q<p+1$
and $D^\pm_\tau(q) = \infty$ for these values. Because of
(\ref{ups10b}) the same happens for $\Delta^{\pm}_\tau(q)$.

Finally, since the density of $\mu$ over $I_j$ grows as $j$, which
is the same as $x^{\frac{1}{p}}$, standard theory gives the
formula for the generalized dimensions $D_\mu(q)$.
 \qed

Notice that divergence of certain moments of entrance and return
times in cylinders for the Manneville --  Pomeau map has been
exhibited in \cite{ugacha}. Further details on this map, and a
{\em local} analysis of return and entrance times can be found in
\cite{xav} in the case when the invariant measure is infinite.


\begin{thebibliography}{99}



\bibitem{miguel}
M. Abadi, {\em Sharp error terms and necessary conditions for
exponential hitting times in mixing processes}, {\em Ann. Probab.}
{\bf  32} (2004) 243--264.
\bibitem{remo2}
R. Badii, A. Politi, {\em Renyi dimensions from local expansion
rates},  {\em Phys. Rev. A} {\bf 35} (1987) 1288--1293.

\bibitem{remo1}
R. Badii, A. Politi, {\em Statistical description of chaotic
attractors: the dimension function},  {\em  J. Stat. Phys.} {\bf
40} (1985)  725--750.

\bibitem{tre} J.M. Barbaroux, F. Germinet, S. Tcheremchantsev,
{\em Generalized fractal dimensions: equivalence and basic
properties}, {\em J. Math. Pures Appl.} {\bf 80} (2001) 977-1012.

\bibitem{xav} X. Bressaud and R. Zweim\"uller,
{\em Non exponential law of entrance times in asymptotically rare
evens for intermittent maps with infinite invariant measure}, {\em
Ann. H. Poincar\'e} {\bf 2} (2001) 501--512.

\bibitem{ugacha}
J.-R. Chazottes and E. Ugalde, {\em Entropy estimation and
fluctuations of hitting and recurrence times for Gibbsian
sources}, {\em Discrete Contin. Dyn. Syst.}  {\bf 5} (2005)
565--586.

\bibitem{colgas} P. Collet, A. Galves and B. Schmitt, {\em Fluctuations of repetition times
for Gibbsian sources}, {\em Nonlinearity} {\bf 12} (1999)
1225--1237.

\bibitem{cutler}
C. D. Cutler, {\em Some results on the behaviour and estimation of
the fractal dimensions of distributions on attractors}, {\em J.
Stat. Phys.} {\bf 62} (1991) 651--708.

\bibitem{gasp} P. Gaspard and X.-J. Wang,
{\em Sporadicity: Between periodic and chaotic dynamical
behaviors}, {\em Proc. Natl. Acad. Sci. U.S.A.} {\bf 85} (1988)
4591--4595.

\bibitem{due} F. Germinet, S. Tcheremchantsev,
{\em Generalized fractal dimensions on the negative axis for
compactly supported measures}, {\em Mathematische Nachrichten}
{\bf 279} (2006) 543--570.

\bibitem{grass}
P. Grassberger, {\em Generalized dimension of strange attractors},
{\em Phys. Lett. A}  {\bf 97}  (1983) 227-–230.

\bibitem{grassp}
P. Grassberger and I. Procaccia,  {\em Characterization of strange
sets},  {\em Phys. Rev. Lett.} {\bf 50} (1983) 346–-349.

\bibitem{nat-0} S. Gratrix and J. N. Elgin, {\em Pointwise dimension of
the Lorenz Attractor}, {\em Phys. Rev. Lett.} {\bf 92} (2004)
14101.

\bibitem{guy}
M. Guysinsky, S. Yaskolko, {\em Coincidence of various dimensions
associated with metrics and measures on metric spaces}, {\em
Discrete Contin. Dynam. Systems} {\bf 3}  (1997) 591--603.

\bibitem{nature} T. C. Halsey and M. H. Jensen, {\em Hurricanes and
Butterflies}, {\em Nature} {\bf 428} (2004) 127.

\bibitem{multif1}
T. C. Halsey, M. Jensen, L. Kadanoff, I. Procaccia, and B.
Shraiman, {\em Fractal measures and their singularities: The
characterization of strange sets}, {\em Phys. Rev. A} {\bf 33}
(1986) 1141--1151.

\bibitem{prl} N. Haydn, J. Luevano, G. Mantica, S. Vaienti, {\em
Multifractal properties of return time statistics}, {\em Phys.
Rev. Lett.} {\bf 88} (2003) 224502.

\bibitem{nico} N. Haydn and S. Vaienti, {\em The limiting distribution and
error terms for return times of dynamical systems},  {\em Discrete
Contin. Dyn. Sys.} {\bf 10} (2004) 589--616.

\bibitem{hents} H. G. E. Hentschel and I. Procaccia, {\em The infinite number of
generalized dimensions of fractals and strange attractors}, {\em
Physica D} {\bf 8} (1983) 435--444.

\bibitem{beno2} M. Hirata, B. Saussol, and S. Vaienti,
{\em Statistics of return times: a general framework and new
applications}, {\em Comm. Math. Phys.} {\bf 206} (1999) 33--55.

\bibitem{kada} M.H. Jensen, L.P. Kadanoff, A. Libchaber, I.
Procaccia and J. Stavans, {\em Global universality at the onset of
chaos: results of a forced Rayleigh--B\'enard experiment}, {\em
Phys. Rev. Lett.} {\bf 55} (1985) 2798--2801.

\bibitem{kac}
M. Kac, ``Probability and related topics in physical sciences'',
AMS, (1959).

\bibitem{sandold} J. Luevano, V. Penne, S. Vaienti, {\em
Multifractal spectrum via return times} Preprint mp-arc 00-232
(2000).

\bibitem{olsen} L. Olsen, {\em First return times: multifractal
spectra and divergence points},  {\em Discrete Contin. Dyn. Sys.}
{\bf 10} (2004) 635--656.

\bibitem{ornstein}
D. S. Ornstein and B. Weiss, {\em Entropy and data compression
schemes}, {\em IEEE Trans. Inform. Theory} {\bf 39} (1993) 78-–83.

\bibitem{sator}
R. Pastor-Satorras and R. H. Riedi, {\em Numerical estimates of
the generalized dimensions of the H\'enon attractor for negative
q}, {\em J. Phys. A: Math. Gen.} {\bf 29} (1996) L391--L398.

\bibitem{pesin-book}
Ya. Pesin, ``Dimension theory in dynamical systems'', The
University of Chicago Press  (1997).

\bibitem{pesjstat}
Ya. Pesin, H. Weiss, {\em A multifractal analysis of equilibrium
measures for conformal expanding maps and Moran-like geometric
constructions}, {\em  J. Stat. Phys.} {\bf  86}  (1997), 233--275.

\bibitem{peters}
K. E. Petersen, ``Ergodic Theory'', Cambridge University Press
(1989).

\bibitem{pmv} Y. Pomeau and P. Manneville, {\em Intermittent transition to
turbulence in dissipative dynamical systems}, {\em Commun. Math.
Phys.} {\bf 74} (1980) 189--197.

\bibitem{rie}
R. H. Riedi, {\em An Improved Multifractal Formalism and
Self-Similar Measures}, {\em  J. Math. Anal. Appl.} {\bf 189},
(1995) 462--490.



\bibitem{beno} B. Saussol, {\em On fluctuations and the exponential statistics of return
times}, {\em Nonlinearity} {\bf 14} (2001) 179--191.

\bibitem{belli} H. Schulz-Baldez and J.
Bellissard, {\em Anomalous transport: a mathematical framework},
{\em Rev. Math. Phys.} {\bf 10} (1998) 1--46.


\bibitem{sandro}
S. Vaienti, {\em  Generalised spectra for the dimensions of
strange sets}, {\em J. Phys. A} {\bf 21}  (1988) 2313--2320.

\bibitem{von} J. von Neumann,
{\em Z\"ur Operatorenmethode in der klassichen Mechanik}, {\em
Ann. of Math.} {\bf 33} (1932) 587--642. (Our map is presented on
page 630).




\end{thebibliography}
\end{document}